\documentclass{amsproc}

\usepackage[english]{babel}
\usepackage[cp1250]{inputenc}
\usepackage{amsmath,amssymb,amsfonts}
\usepackage{url}
\usepackage[dvipsnames,usenames]{color}
\usepackage{array}
\usepackage{accents}
\usepackage{tabu}
\usepackage{mathdots}

\usepackage{tikz}
\usetikzlibrary{arrows,positioning,matrix} 
\tikzset{
    >=stealth',
    punkt/.style={
           rectangle,
           rounded corners,
           draw=black, thick,
           text width=6.5em,
           minimum height=2em,
           text centered},
    punkt2/.style={
           rectangle,
           rounded corners,
           draw=black, thick,
           text width=8.5em,
           minimum height=2em,
           text centered},
    pil/.style={
           ->,
           thick,
           shorten <=2pt,
           shorten >=2pt,}
}

\newcolumntype{L}{>{\displaystyle}l}
\newcolumntype{C}{>{\displaystyle}c}
\newcolumntype{R}{>{\displaystyle}r}

\renewcommand{\tfrac}{\genfrac{}{}{}1}

\newcommand{\R}{\mathbb R}
\newcommand{\N}{\mathbb N}

\newcommand{\C}{\mathbb C}

\renewcommand{\Im}{\mathrm{Im}}

\def\til{\widetilde}
\newcommand{\Res}{\mathrm{Res}}
\def\I{\mathfrak{i}}
\newcommand{\diff}{\mathrm{d}}
\renewcommand{\bar}{\overline}

\newcommand{\ie}{\textit{i.e.}\/ } 
 
\newcommand{\cf}{\textit{cf.}\/ }

\renewcommand{\vec}[1]{\accentset{\rightharpoonup}{#1}} 
\newcommand{\mat}[1]{#1} 
\newcommand{\tilvec}[1]{\til{#1}} 
\newcommand{\hatvec}[1]{\hat{#1}}

\allowdisplaybreaks

\theoremstyle{definition}
\newtheorem{define}{Definition}[section]
\newtheorem{example}[define]{Example}
\newtheorem{remark}[define]{Remark}

\theoremstyle{plain}

\newtheorem{thm}[define]{Theorem}
\newtheorem{prop}[define]{Proposition}
\newtheorem{coro}[define]{Corollary}

\numberwithin{equation}{section}

\begin{document}

\title[The convex combination problem]{A subclass of boundary measures and the convex combination problem for Herglotz-Nevanlinna functions in several variables}

\author{Mitja Nedic}
\address{Mitja Nedic, Department of Mathematics, Stockholm University, SE-106 91 Stockholm, Sweden, orc-id: 0000-0001-7867-5874}
\curraddr{}
\email{mitja@math.su.se}
\thanks{The author is supported by the Swedish Foundation for Strategic Research, grant nr. AM13-0011.}


\date{17-03-2017}

\begin{abstract}
In this paper, we begin by investigating a particular subclass of boundary measures of Herglotz-Nevanlinna functions in two variables. Based on this, we then proceed to solve the convex combination problem for Herglotz-Nevanlinna functions in several variables.
\end{abstract}

\maketitle

\section{Introduction}\label{sec:introduction}

Herglotz-Nevanlinna functions are holomorphic functions defined on the poly-upper half-plane having non-negative imaginary part. In the classical case of one complex variable, these functions have proven to be most useful, both within the fields of mathematical analysis and electromagnetic engineering. This development started around 100 years ago with Rolf H. Nevanlinna's work on the Stieltjes moment problem \cite{Nevanlinna1922}. Since then, these functions have found their home, among others places, in spectral theory \cite{KacKrein1974,Simon1998}, the moment problem \cite{Akhiezer1965,Nevanlinna1922,Simon1998} and convex optimization \cite{Ivanenko+etal2017} on the mathematical side, as well as  derivation of physical bounds \cite{Bernland2011,Milton2016,Milton2002}, homogenization of two-component media \cite{GoldenPapanicolaou1983,Milton2002} and circuit synthesis \cite{Brune1931} on the engineering side. This has primary been possible due to the powerful integral representation theorem for these class of functions \cite{Cauer1932}, \cf Theorem \ref{thm:intRep_1var}.

The class of Herglotz-Nevanlinna functions in several variables is, however, a slightly newer consideration, appearing first in the works of Vladimirov and his collaborators \cite{Vladimirov1974,Vladimirov1979,Vladimirov1969} in the 1970s. This class of functions has, so far, proven a bit less prominent in both mathematics and applications, but relates nonetheless to multidimensional passive systems \cite{Vladimirov1969} and homogenization of multicomponent media \cite{GoldenPapanicolaou1985,Milton2016,Milton2002}. However, the class of Herglotz-Nevanlinna functions in several variables has seen a renewed interest in the last few years, especially from a pure mathematical perspective, with results concerning both integral representations of this class of functions \cite{LugerNedic2016,LugerNedic2017}, as well as operator representations \cite{agler2,agler}.

Consider now the following problem. Suppose we are given a Herglotz-Nevan\-linna function $q$ in one variable, which we use to build a new Herglotz-Nevanlinna function $\til{q}$ in several variables by replacing the argument of the function $q$ with a convex combination of several variables, \ie
$$\til{q}\colon (z_1,z_2,\ldots,z_n) \mapsto q(k_1z_1+k_2z_2+\ldots+ k_nz_n),$$
where the coefficients $k_j$ describe the convex combination. We are then interested in relating the parameters of the integral representation of the function $q$ to the parameters of the integral representation of the function $\til{q}$ in the most explicit way possible. We refer to this conundrum as the \emph{convex combination problem}, with the \emph{arithmetic mean problem} being the obvious special case of the above question. In the beginning, only the arithmetic mean problem was considered and an answer obtained, but questions arising from a seminary discussion quickly encouraged the consideration of the more general problem.

The motivation behind these types of problems comes mainly from two viewpoints. On one side, there is the desire to have a table of explicit integration formulas, not unlike those that one may find in Abramowitz and Stegun's classical work \cite{AbramowitzStegun1964}, or, perhaps, in King's encyclopedia on the Hilbert transform \cite{King2009}. On the other side, we wish to relate the data of one Herglotz-Nevanlinna function to another in the case when these functions are related by some identity. This is elaborated upon later in Section \ref{sec:convex_solution}.

In this paper, we provide a completely explicit answer to the complex combination problem in full generality using the boundless power of classic residue calculus in one complex variable. This is presented in Theorem \ref{thm:Nvar_convex}, which is, therefore, the main result of this paper. 

After the introduction in Section \ref{sec:introduction}, we continue with a short review of the integral representation formula for Herglotz-Nevanlinna function, presented in Section \ref{sec:int_representations}. In Section \ref{sec:subclass}, we consider a particular subclass of boundary measures of Herglotz-Nevanlinna functions in two variables, that turns out to be the starting point of the solution of the arithmetic mean problem in two variables. This solution is presented in Section \ref{sec:convex_solution}, along with the solutions of the convex combination problem in two variables and the solutions of both the arithmetic mean and convex combination problems in full generality.

\section{Integral representations of Herglotz-Nevanlinna functions}\label{sec:int_representations}

Let us begin by recalling some know results about integral representations of Herglotz-Nevanlinna functions. Throughout this paper, we will denote by the letter $z$ variables which lie in the upper half-plane, while the letter $t$ is reserved for real-valued variables. Recall also that the poly-upper half-plane is defined by
$$\C ^{+n} := (\C^+)^n = \{\vec{z} \in \C^n~|~\Im[z_j] > 0~\text{for all}~j = 1,\ldots,n\}.$$

The integration kernel that we will be considering is the kernel $K_n$, defined for $\vec{z} \in \C^{+n}$ and $\vec{t} \in \R^n$ as
\begin{equation}\label{eq:kernel_Kn}
K_n(\vec{z},\vec{t}) := \I\left(\frac{2}{(2\I)^n}\prod_{j=1}^{n}\left(\frac{1}{t_j-z_j}-\frac{1}{t_j+\I}\right) - \frac{1}{(2\I)^n}\prod_{j=1}^{n}\left(\frac{1}{t_j-\I}-\frac{1}{t_j+\I}\right)\right),
\end{equation}
and can be equivalently written as
\begin{equation}\label{eq:kernel_Kn_v2}
K_n(\vec{z},\vec{t}) = \frac{\I^{3n+1}\prod_{j=1}^n(t_j-\I)(z_j+\I)-2^{n-1}\I\prod_{j=1}^n(t_j-z_j)}{2^{n-1}\prod_{j=1}^n(t_j-z_j)(t_j-\I)(t_j+\I)}.
\end{equation}
The class of all Herglotz-Nevanlinna function can then be completely characterized via an integral representation formula \cite[Theorem 4.1]{LugerNedic2017}, as described by the following theorem. 

\begin{thm}\label{thm:intRep_Nvar}
A function $q:\C^{+n} \to \C$ is a Herglotz-Nevanlinna function if and only if it admits, for $\vec{z} \in \C^{+n}$, a representation of the form
\begin{equation}\label{eq:intRep_Nvar}
q(\vec{z}) = a + \sum_{\ell=1}^nb_\ell\:z_\ell + \frac{1}{\pi^n}\int_{\R^n}K_n(\vec{z},\vec{t})\diff\mu(\vec{t}),
\end{equation}
where $a \in \R$, $b_\ell \geq 0$ for all $\ell = 1,\ldots,n$ and $\mu$ is a positive Borel measure on $\R^n$ satisfying the growth condition
\begin{equation}\label{eq:growth_Nvar}
\int_{\R^n}\prod_{j=1}^n\frac{1}{1+t_j^2}\diff\mu(\vec{t}) < \infty
\end{equation}
and the Nevanlinna condition, \ie
\begin{equation}\label{eq:Nevan_Nvar}
\sum_{\substack{\vec{\rho} \in \{-1,0,1\}^n \\ -1\in\vec{\rho} \wedge 1\in\vec{\rho}}}\int_{\R^n}N_{\rho_1,1}N_{\rho_2,2}\ldots N_{\rho_n,n}\diff\mu(\vec{t}) = 0
\end{equation} 
for all $\vec{z} \in \C^{+n}$, where the factors $N_{k,j}$ are defined as
\begin{equation*}
\label{eq:Nevan_factors}
N_{-1,j} := \frac{1}{t_j - z_j} - \frac{1}{t_j - \I},~N_{0,j} := \frac{1}{t_j - \I} - \frac{1}{t_j + \I},~N_{1,j} :=\frac{1}{t_j + \I} - \frac{1}{t_j - \bar{z_j}}.
\end{equation*}
\end{thm}

\begin{remark}\label{rem:intRep_unique}
It can be shown that the above correspondence is, in fact, a bijection. That is to say, the parameters $a$, $\vec{b}$ and $\mu$ are unique for a given Herglotz-Nevanlinna function $q$, and conversely, a different choice of parameters corresponds to a different function \cite{LugerNedic2017}. Therefore, for simplicity, we often say that a function $q$ is represented by the \emph{data} $(a,\vec{b},\mu)$. This can even be improved in the following way. Suppose that $a_1 \in \R, \vec{b}_1 \in \R^n$ and a positive Borel measure $\mu_1$ on $\R^n$ are such that they give a Herglotz-Nevanlinna function when plugged into representation \eqref{eq:intRep_Nvar}. Then, they must satisfy the conditions of Theorem \ref{thm:intRep_Nvar} and are, in fact, equal to the data of the function in question \cite[Corolarry 4.7]{LugerNedic2017}.
\end{remark}

\begin{remark}
It can be shown that the measure $\mu$ from Theorem \ref{thm:intRep_Nvar} is in fact the limit of the function $\Im[q]$ as we approach $\R^n$ from $\C^{+n}$. Therefore, the measure $\mu$ is also called both the \emph{representing measure} and the \emph{boundary measure} of a function $q$ \cite{LugerNedic2016,LugerNedic2017}.
\end{remark}

When $n=2$, we note that the growth condition \eqref{eq:growth_Nvar} becomes
\begin{equation}\label{eq:growth_2var}
\int_{\R^2}\frac{1}{(1+t_1^2)(1+t_2^2)}\diff\mu(\vec{t}) < \infty
\end{equation}
and that the Nevanlinna condition \ref{eq:Nevan_Nvar} is then equivalent to the condition that
\begin{equation}\label{eq:Nevan_2var}
\int_{\R^2}\frac{1}{(t_1-z_1)^2(t_2-\bar{z_2})^2}\diff\mu(\vec{t}) = 0
\end{equation}
for any $(z_1,z_2) \in \C^{+2}$ \cite[Theorem 5.1]{LugerNedic2017}. It is because of this equivalence that condition \eqref{eq:Nevan_2var} is, for simplicity, also referred to as the Nevanlinna condition (in two variables).

We note also that in the case $n=1$, Theorem \ref{thm:intRep_Nvar} reduces the classical theorem attributed to Nevanlinna, presented in its current form by Cauer \cite{Cauer1932}.

\begin{thm}\label{thm:intRep_1var}
A function $q:\C^{+} \to \C$ is a Herglotz-Nevanlinna function in variables if and only if it admits, for $z \in \C^+$, a representation of the form
\begin{equation}\label{eq:intRep_1var}
q(z) = a + bz + \frac{1}{\pi}\int_{\R}K_1(z,t)\diff\mu(t),
\end{equation}
where $a \in \R$, $b \geq 0$ and $\mu$ is a positive Borel measure on $\R$ satisfying the growth condition
\begin{equation}\label{eq:growth_1var}
\int_{\R}\frac{1}{1+t^2}\diff\mu(t) < \infty.
\end{equation}
\end{thm}

\section{A special class of boundary measures for functions of two variables}\label{sec:subclass}

Let us consider now a particular subclass of boundary measures of Herglotz-Nevanlinna functions in two variables. The introduction of our particular subclass of measures is motivated by the following example.

\begin{example}\label{ex:motivation_2var}
Let $q$ be a Herglotz-Nevanlinna function in one variable, represented by the data $(a,b,\mu)$. Consider now two Herglotz-Nevanlinna functions $\til{q}_1,\til{q}_2$ in two variables, defined by
$$\til{q}_1\colon (z_1,z_2) \mapsto q(z_1)$$
and
$$\til{q}_2\colon (z_1,z_2) \mapsto q(z_2),$$
respectively. It can be shown that the function $\til{q}_1$ is represented by the data $(a,(b,0),\mu \otimes \lambda_\R)$, while the function $\til{q}_2$ is represented by the data $(a,(0,b),\lambda_\R \otimes \mu)$. Here, $\lambda_\R$ denotes the Lebesgue measure on $\R$. This follows from the fact that integrating the kernel $K_n$ once with respect to $\diff t_j$ gives a constant multiple of $K_{n-1}$ with the $j$-th variable missing \cite[Example 3.4]{LugerNedic2017}.\hfill$\lozenge$
\end{example}

Given the above example, we are led to conjecture that a Herglotz-Nevanlinna function given by $(z_1,z_2) \mapsto q(k_1 z_1 + k_2 z_2)$ with $k_1,k_2 > 0, k_1+k_2=1$, should have a boundary measure that is "somewhere in between" $\mu \otimes \lambda_\R$ and $\lambda_\R \otimes \mu$. We formalize this idea by introducing the following class of Borel measures on $\R^2$.

First, let $\alpha,\beta,\gamma,\delta \in \R$ and let $\mu_1$ be a positive Borel measure on $\R$. We then consider the Borel measure $\mu$ on $\R^2$, defined for any Borel measurable subset $U \subseteq \R^2$ as
\begin{equation}\label{eq:measure_2var_general}
\mu(U) := \int_{\R}\left(\int_{\R}\chi_{U}(\alpha t_1 + \beta t_2,\gamma t_1 + \delta t_2)\diff t_2\right)\diff\mu_1(t_1).
\end{equation}

\begin{remark}\label{rem:integration_order}
Throughout this paper, we never discuss what happens if the order of integration in the above definition of the measure $\mu$ is reversed. For us, the inner integral in formula \eqref{eq:measure_2var_general} is always taken first and is always with respect to the Lebesgue measure.
\end{remark}

We now ask the question whether measures on $\R^2$ of the type \eqref{eq:measure_2var_general} can be representing measures of Herglotz-Nevanlinna functions in two variables. As is turns out, the answer can be summarized by the following theorem.

\begin{thm}\label{thm:measure_2var}
A positive Borel measure $\mu$ of the form \eqref{eq:measure_2var_general} is the representing measure of some Herglotz-Nevanlinna function in two variables if and only if one of the following cases holds:
\begin{itemize}
\item[(i.1)]{$\alpha = 0$, $\beta = 0$, $\delta \neq 0$ and $\mu_1$ is a finite positive Borel measure on $\R$,}
\item[(i.2)]{$\alpha \neq 0$, $\beta = 0$, $\delta \neq 0$ and $\mu_1$ is a positive Borel measure on $\R$ satisfying the growth condition \eqref{eq:growth_1var},}
\item[(ii.1)]{$\beta \neq 0$, $\gamma = 0$, $\delta = 0$ and $\mu_1$ is a finite positive Borel measure on $\R$,}
\item[(ii.2)]{$\beta \neq 0$, $\gamma \neq 0$, $\delta = 0$ and $\mu_1$ is a positive Borel measure on $\R$ satisfying the growth condition \eqref{eq:growth_1var},}
\item[(iii.1.a)]{$\beta\delta < 0$, $\alpha\delta-\beta\gamma = 0$ and $\mu_1$ is a finite positive Borel measure on $\R$,}
\item[(iii.1.b)]{$\beta\delta < 0$, $\alpha\delta-\beta\gamma \neq 0$ and $\mu_1$ is a positive Borel measure on $\R$ satisfying the growth condition \eqref{eq:growth_1var},}
\item[(iii.2.a)]{$\beta\delta > 0$, $\alpha\delta-\beta\gamma = 0$ and $\mu_1$ is identically zero,}
\item[(iii.2.b)]{$\beta\delta > 0$, $\alpha\delta-\beta\gamma \neq 0$ and $\mu_1$ is a positive Borel measure on $\R$ satisfying the growth condition \eqref{eq:growth_1var} and the condition that
\begin{equation*}
\int_\R\frac{1}{((\alpha\delta-\beta\gamma)t_1-\delta z_1 + \beta \bar{z_2})^3}\diff\mu_1(t_1) = 0
\end{equation*}
for all $z_1,z_2 \in \C^+$.}
\end{itemize}
\end{thm}

The remainder of this section is devoted to the proof of this theorem and is divided into two propositions. Proposition \ref{prop:measure_growth} first characterizes which measures of the form \eqref{eq:measure_2var_general} satisfy the growth condition \eqref{eq:growth_2var}, while Proposition \ref{prop:measure_nevan} then characterizes which measures of the form \eqref{eq:measure_2var_general} satisfy the Nevanlinna condition condition \eqref{eq:Nevan_2var}. Combining these results gives Theorem \ref{thm:measure_2var}.

We begin now with the first of the aforementioned propositions.

\begin{prop}\label{prop:measure_growth}
A measure $\mu$ of the type \eqref{eq:measure_2var_general} satisfies the growth condition \eqref{eq:growth_2var} if and only one of the following cases holds:
\begin{itemize}
\item[(i.1)]{$\alpha = 0$, $\beta = 0$, $\delta \neq 0$ and $\mu_1$ is a finite positive Borel measure on $\R$,}
\item[(i.2)]{$\alpha \neq 0$, $\beta = 0$, $\delta \neq 0$ and $\mu_1$ is a positive Borel measure on $\R$ satisfying the growth condition \eqref{eq:growth_1var},}
\item[(ii.1)]{$\beta \neq 0$, $\gamma = 0$, $\delta = 0$ and $\mu_1$ is a finite positiveBorel measure on $\R$,}
\item[(ii.2)]{$\beta \neq 0$, $\gamma \neq 0$, $\delta = 0$ and $\mu_1$ is a positive Borel measure on $\R$ satisfying the growth condition \eqref{eq:growth_1var},}
\item[(iii.1.a)]{$\beta\delta < 0$, $\alpha\delta-\beta\gamma = 0$ and $\mu_1$ is a finite positive Borel measure on $\R$,}
\item[(iii.1.b)]{$\beta\delta < 0$, $\alpha\delta-\beta\gamma \neq 0$ and $\mu_1$ is a positive Borel measure on $\R$ satisfying the growth condition \eqref{eq:growth_1var},}
\item[(iii.2.a)]{$\beta\delta > 0$, $\alpha\delta-\beta\gamma = 0$ and $\mu_1$ is a finite positive Borel measure on $\R$,}
\item[(iii.2.b)]{$\beta\delta > 0$, $\alpha\delta-\beta\gamma \neq 0$ and $\mu_1$ is a positive Borel measure on $\R$ satisfying the growth condition \eqref{eq:growth_1var}.}
\end{itemize}
\end{prop}

\proof
The total integral, appearing in the growth condition \eqref{eq:growth_2var} for a measure $\mu$ of the form \eqref{eq:measure_2var_general}, is equal to
\begin{multline}\label{eq:growth_total_int}
\int_{\R^2}\frac{1}{(1+t_1^2)(1+t_2^2)}\diff\mu(\vec{t}) \\
= \int_\R\left(\int_\R\frac{1}{(\alpha t_1+\beta t_2-\I)(\alpha t_1+\beta t_2+\I)(\gamma t_1+\delta t_2-\I)(\gamma t_1+\delta t_2+\I)}\diff t_2\right)\diff\mu_1(t_1).
\end{multline}
We now investigate the finiteness of this integral with respect to the numbers $\beta$ and $\delta$.

Observe first that the inner integral
\begin{equation}\label{eq:growth_inner_int}
\int_\R\frac{1}{(\alpha t_1+\beta t_2-\I)(\alpha t_1+\beta t_2+\I)(\gamma t_1+\delta t_2-\I)(\gamma t_1+\delta t_2+\I)}\diff t_2
\end{equation}
cannot be finite unless at least one of the numbers $\beta$ and $\delta$ is non-zero. If one of the numbers $\beta$ and $\delta$ is equal to zero, then the integral \eqref{eq:growth_inner_int} is, in the case $\beta = 0$ and $\delta \neq 0$, equal to
$$\frac{1}{\alpha^2t_1^2+1}\int_\R\frac{1}{(\gamma t_1 + \delta t_2)^2 + 1}\diff t_2 = \frac{1}{\alpha^2t_1^2+1} \cdot \frac{\pi}{|\delta|}.$$
The case $\beta \neq 0$ and $\delta = 0$ is treated analogously. Therefore, when one of the numbers $\beta$ and $\delta$ is non-zero, the total integral \eqref{eq:growth_total_int} becomes finite if and only if one the first four cases happens.

If both numbers $\beta$ and $\delta$ are non-zero, we are left to consider the cases $\beta\delta < 0$ and $\beta\delta > 0$. We begin by investigating the case $\beta\delta < 0$ by using standard residue calculus to calculate the inner integral \eqref{eq:growth_inner_int}. Let now 
$$F(\tau) := \frac{1}{(\alpha t_1+\beta \tau-\I)(\alpha t_1+\beta \tau+\I)(\gamma t_1+\delta \tau-\I)(\gamma t_1+\delta \tau+\I)}$$
be an auxiliary function, where the parameters $\alpha,\beta,\gamma,\delta \in \R$, as well as on $t_1 \in \R$, are fixed. We note now that the integral
$$\int_\R F(\tau)\diff \tau$$
is well-defined since since the function $F$ is a rational function with a constant numerator, while the denominator is a polynomial of degree $4$. Note that this observation is valid independently of the particular values of the parameters $\alpha,\beta,\gamma,\delta$ and $t_1$. Next, observe that the function $F$ has singularities at the points
$$\frac{\I-\alpha t_1}{\beta},\frac{-\I-\alpha t_1}{\beta},\frac{\I-\gamma t_1}{\delta},\frac{-\I-\gamma t_1}{\delta} \in \C\setminus\R.$$

Consider now the case when $\beta > 0$ and $\delta < 0$ and take
$$R > \max\left\{\frac{|\I-\alpha t_1|}{\beta},\frac{|\I+\gamma t_1|}{-\delta}\right\}.$$
Let also $\Gamma_R^+$ be the standard upper half-circle contour in $\C$, \ie the curve consisting of the interval $[-R,R]$ and the curve $\gamma_R^+$, which is the upper half-circle of radius $R$ centered at $0$ (note that the curve $\gamma_R^+$ has no connection to the number $\gamma$), oriented counter-clockwise. Then, due to the rational form of the function $F$, it likewise holds that
$$\lim\limits_{R \to \infty}\int_{\Gamma_R^+}F(\tau)\diff\tau = \lim\limits_{R \to \infty}\left(\int_{-R}^R+\int_{\gamma_R^+}\right)F(\tau)\diff\tau = \int_\R F(\tau)\diff\tau,$$
while, by the residue theorem, it holds that
$$\begin{array}{RCL}
\multicolumn{3}{L}{\int_{\Gamma_R^+}F(\tau)\diff\tau = 2\pi\I\big(\Res(F;\tfrac{\I-\alpha t_1}{\beta})+\Res(F;\tfrac{-\I-\gamma t_1}{\delta})\big)} \\[0.5cm]
~ & = & 2\pi\I\bigg(\lim\limits_{\tau \to \tfrac{\I-\alpha t_1}{\beta}}F(\tau)(\tau - \tfrac{\I-\alpha t_1}{\beta}) + \lim\limits_{\tau \to \tfrac{-\I-\gamma t_1}{\delta}}F(\tau)(\tau - \tfrac{-\I-\gamma t_1}{\delta})\bigg) \\[0.5cm]
~ & = & 2\pi\I\bigg(\frac{\beta}{2\I}\:\frac{1}{(t_1(\beta\gamma-\alpha\delta)-\I(\beta-\delta))(t_1(\beta\gamma-\alpha\delta)+\I(\beta+\delta))} \\[0.5cm]
~ & ~ & - \frac{\delta}{2\I}\:\frac{1}{(t_1(\beta\gamma-\alpha\delta)+\I(\beta+\delta))(t_1(\beta\gamma-\alpha\delta)+\I(\beta-\delta))}\bigg) \\[0.5cm]
~ & = & \frac{\pi(\beta-\delta)}{t_1^2(\beta\gamma-\alpha\delta)^2+(\beta-\delta)^2}.
\end{array}$$
Observe here that $\beta-\delta \neq 0$ since we are in working with the case when $\beta > 0$ and $\delta < 0$. Thus, the total integral \eqref{eq:growth_total_int} becomes
$$\int_\R\frac{\pi(\beta-\delta)}{t_1^2(\beta\gamma-\alpha\delta)^2+(\beta-\delta)^2}\diff\mu_1(t_1)$$
and is finite if and only if one of the cases (iii.1.a) or (iii.1.b) happens. The case $\beta < 0$ and $\delta > 0$ is considered analogously.

Finally, we see that, in the case $\beta\delta > 0$, the total integral \eqref{eq:growth_total_int} is finite if one of the cases (iii.2.a) or (iii.2.b) happens through and analogous application of the residue theorem.
\endproof

\begin{prop}\label{prop:measure_nevan}
A measure $\mu$ of the type \eqref{eq:measure_2var_general} satisfies the Nevanlinna condition \eqref{eq:Nevan_2var} if and only if one of the following cases holds:
\begin{itemize}
\item[(i)]{$\beta = 0$, $\delta \neq 0$ and $\mu_1$ is a positive Borel measure on $\R$,}
\item[(ii)]{$\beta \neq 0$, $\delta = 0$ and $\mu_1$ is a positive Borel measure on $\R$,}
\item[(iii.1)]{$\beta\delta < 0$ and $\mu_1$ is a positive Borel measure on $\R$,}
\item[(iii.2.a)]{$\beta\delta > 0$, $\alpha\delta-\beta\gamma = 0$ and $\mu_1$ is identically zero,}
\item[(iii.2.b)]{$\beta\delta > 0$, $\alpha\delta-\beta\gamma \neq 0$ and $\mu_1$ is a positive Borel measure on $\R$ satisfying the condition that
\begin{equation}\label{eq:measure_nevan_case}
\int_\R\frac{1}{((\alpha\delta-\beta\gamma)t_1-\delta z_1 + \beta \bar{z_2})^3}\diff\mu_1(t_1) = 0
\end{equation}
for all $z_1,z_2 \in \C^+$.
}
\end{itemize}
\end{prop}

\proof
The total integral, appearing in the Nevanlinna condition \eqref{eq:Nevan_2var} for a measure $\mu$ of the form \eqref{eq:measure_2var_general}, is equal to
\begin{multline}\label{eq:nevan_total_int}
\int_{\R^2}\frac{1}{(t_1-z_1)^2(t_2-\bar{z_2})^2}\diff\mu(\vec{t}) \\
= \int_{\R}\left(\int_{\R}\frac{1}{(\alpha t_1+\beta t_2-z_1)^2(\gamma t_1+\delta t_2-\bar{z_2})^2}\diff t_2\right)\diff\mu_1(t_1)
\end{multline}
Similarly to the previous proof, we now investigate when this integral is identically equal to zero with respect to the numbers $\beta$ and $\delta$.

We observe first that the inner integral
\begin{equation}\label{eq:nevan_inner_int}
\int_{\R}\frac{1}{(\alpha t_1+\beta t_2-z_1)^2(\gamma t_1+\delta t_2-\bar{z_2})^2}\diff t_2
\end{equation}
cannot be finite unless at least one of the numbers $\beta$ and $\delta$ is non-zero.

If one of the numbers $\beta$ and $\delta$ is equal to zero, then the integral \eqref{eq:nevan_inner_int} becomes trivial to compute with the help of the residue theorem. For example, in the case $\beta = 0$ and $\delta \neq 0$, the inner integral \eqref{eq:nevan_inner_int} becomes
$$\frac{1}{(\alpha t_1 - z_1)^2}\int_\R\frac{1}{(\gamma t_1 + \delta t_2 - \bar{z_2})^2}\diff t_2.$$
Its integrand is, with respect to the variable $t_2$, a rational function whose denominator is a polynomial of degree at $2$. This allows for the use of the residue theorem. Since the integrand has only one singularity in the complexified $t_2$-variable, we see quickly that the inner integral \eqref{eq:nevan_inner_int} is identically zero in this case. The total integral \eqref{eq:nevan_total_int} is therefore also identically zero for any positive Borel measure $\mu_1$. The case $\beta \neq 0$, and $\delta = 0$ can be considered completely analogously. This give the first two cases of the proposition. Here, it is also important to remember that we always abide by Remark \ref{rem:integration_order}.

If now both numbers $\beta$ and $\delta$ are non-zero, we are left to consider the cases $\beta\delta < 0$ and $\beta\delta > 0$. We begin by investigating the case $\beta\delta < 0$, where we, again, use standard residue calculus to calculate the inner integral \eqref{eq:nevan_inner_int}. To that end, define an auxiliary function $G$ as
$$G(\tau) := \frac{1}{(\alpha t_1+\beta\tau-z_1)^2(\gamma t_1+\delta\tau-\bar{z_2})^2},$$
where the parameters $\alpha,\eta,\gamma,\delta,t_1 \in \R$ are fixed. The function $G$ has singularities at the points
$$\frac{z_1-\alpha t_1}{\beta},\frac{\bar{z_2}-\gamma t_1}{\delta} \in \C\setminus\R.$$
Since $\beta\delta < 0$, then these singularities both lie in the same half-plane. More precisely, if $\beta >0$ and $\delta < 0$ then both lie in the upper half-plane, otherwise they both lie in the lower half-plane. 

We note also that the integral
$$\int_\R G(\tau)\diff \tau$$
is well-defined since the function $G$, similarly to the function $F$ in the previous proof, is a rational function with a constant numerator, while its denominator is a polynomial of degree $4$. Note that this observation is valid independently of the particular values of the parameters $\alpha,\beta,\gamma,\delta$ and $t_1$. Take now
$$R > \max\left\{\frac{|z_1-\alpha t_1|}{\beta},\frac{|\bar{z_2}-\gamma t_1|}{-\delta}\right\}$$
and consider first the case when $\beta > 0$ and $\delta < 0$. Let $\Gamma_R^-$ be the standard lower half-circle contour in $\C$, \ie the curve consisting of the interval $[-R,R]$ and the curve $\gamma_R^-$, which is the lower half-circle of radius $R$ centered at $0$ (note that the curve $\gamma_R^-$ has no connection to the number $\gamma$), oriented clockwise. Then, due to the rational form of the function $G$, it holds that
\begin{equation}\label{eq:nevan_aux_integral}
\lim\limits_{R \to \infty}\int_{\Gamma_R^-}G(\tau)\diff\tau = \lim\limits_{R \to \infty}\left(\int_{-R}^R+\int_{\gamma_R^-}\right)G(\tau)\diff\tau = \int_\R G(\tau)\diff\tau,
\end{equation}
while, by the residue theorem, we conclude that
$$\int_{\Gamma_R^-}G(\tau)\diff\tau = 0.$$
We note that the case $\beta < 0$ and $\delta > 0$ is done completely analogously using the standard upper half-circle contour. Thus, the total integral \eqref{eq:nevan_total_int} is identically zero, in this case, if and only if case (iii.1) happens.

Therefore, it remains to consider the case $\beta\delta > 0$. Using the same auxiliary function $G$ as before, relation \eqref{eq:nevan_aux_integral} still holds. On the other hand, we calculate using the residue theorem that, in the case $\beta > 0$ and $\delta > 0$,
\begin{multline*}
\int_{\Gamma_R^-}G(\tau)\diff\tau = -2\pi\I\Res(G;\tfrac{\bar{z_2}-\gamma t_1}{\delta}) \\
= -2\pi\I\lim\limits_{\tau \to \tfrac{\bar{z_2}-\gamma t_1}{\delta}}(G(\tau)(\tau - \tfrac{\bar{z_2}-\gamma t_1}{\delta})^2)' = \frac{4\pi\I\beta\delta}{((\alpha\delta-\beta\gamma)t_1-\delta z_1 + \beta \bar{z_2})^3}.
\end{multline*}
Here, the accent $'$ denotes the derivative with respect to the $\tau$-variable. The case $\beta < 0$ and $\delta < 0$ can be treated analogously. Thus, the total integral \eqref{eq:nevan_total_int} is identically zero, in this case, if and only if one of the cases (iii.2.a) or (iii.2.b) happens. This finishes the proof.
\endproof

\section{The solution of the convex combination problem}\label{sec:convex_solution}

A common question concerning integral representations of Herglotz-Nevanlinna functions is to relate the data of one function to the data of another, when the  two functions are related by a certain identity. A simple staring example is to consider a Herglotz-Nevanlinna function of one variable, represented by the data $(a,b,\mu)$ and a Herglotz-Nevanlinna function of two variables $\til{q}$, represented by the data $(\til{a},\tilvec{b},\til{\mu})$. Suppose that this functions are related by the identity
$$\til{q}(z_1,z_2) = 1 + 2z_2 + 3q(z_1)$$
for all $(z_1,z_2) \in \C^{+2}$. Then, by writing out both functions using their respective integral representations, one sees that $\til{a} = 1 + 3a$, $\tilvec{b} = (3b,2)$ and $\til\mu = 3\mu \otimes \lambda_\R$. Of course, if the identity relating the functions $q$ and $\til{q}$ is more complicated, for example
$$\til{q}(z_1,z_2) = 1 - \frac{1}{q(z_1 - \frac{1}{z_2} + \I) + \I},$$
it may be utterly impossible to say anything about the relations between the data of the two functions.

Often, we adopt a different point of view to the above problem. In particular, we consider instead the function $\til{q}$ as being built with the help of the function $q$ and ask to relate the data of the starting function to the data of the new function. It is in this spirit that we also view our main problem of interest, namely, the convex combination problem.

We recall from Section \ref{sec:introduction} that the convex combination problem supposes that we are given a Herglotz-Nevnanlinna function $q$ in one variable, which is then used to build a new Herglotz-Nevanlinna function $\til{q}$ in several variables by replacing the argument of the function $q$ with a convex combination of several independent variables, \ie
$$\til{q}\colon (z_1,z_2,\ldots,z_n) \mapsto q(k_1z_1+k_2z_2+\ldots +k_nz_n),$$
where the coefficients $k_\ell > 0$ are such that $k_1+k_2 + \ldots + k_n = 1$. Later, in Corollary \ref{coro:intRep_convex_general}, we will remove the constraint that the coefficients $k_\ell$ are positive.

We are now interested in writing the data $(\til{a},\tilvec{b},\til{\mu})$, corresponding to the function $\til{q}$,  in terms of the data $(a,b,\mu)$, corresponding to the function $q$. The answer to this question is the main result of this paper and is presented by the following theorem.

\begin{thm}\label{thm:Nvar_convex}
Let $q$ be a Herglotz-Nevanlinna function in one variable, represented by the data $(a,b,\mu)$. Let now $n \geq 2$ and $k_\ell > 0$ for $\ell=1,2,\ldots,n$, such that $k_1 + k_2 + \ldots k_n = 1$. Then, the function $\til{q}\colon\C^{+n} \to \C$, defined by
$$\til{q}\colon (z_1,z_2,\ldots,z_n) \mapsto q(k_1z_1 + k_2z_2 + \ldots + k_nz_n),$$
is a Herglotz-Nevanlinna function represented by the data $(\til{a},\til{b},\til{\mu})$, where $\til{a} = a$, $\tilvec{b} = (k_1b,k_2b,\ldots,k_nb)$ and $\til{\mu}$ is a positive Borel measure on $\R^n$, defined for any Borel measurable subset $U \subseteq \R^n$ as
\begin{multline}\label{eq:measure_convex_Nvar}
\til{\mu}(U) \\ := \beta_n\int_\R\left(\int_{\R^{n-1}}\chi_U(t_1 - b_1t_2,\ldots,t_1-b_{n-1}t_n,t_1+t_2+\ldots+t_n)\diff t_n\ldots \diff t_2\right)\diff\mu(t_1).
\end{multline}
Here, the numbers $b_j$, $j = 1,2,\ldots,n-1$, are related to the coefficients $k_\ell$ by relations \eqref{eq:convex_coeff_Nvar} and the number $\beta_n$ is defined by relation \eqref{eq:betaN}.
\end{thm}

We will prove this theorem by showing that the parameters $\til{a},\tilvec{b}$ and $\til{\mu}$, as specified by the theorem, give back the function $\til{q}$ when plugged into the integral representation formula \eqref{eq:intRep_Nvar}, relying also on the uniqueness of the parameters as discussed in Remark \ref{rem:intRep_unique}. This, however, requires substantial calculations involving the use of standard residue calculus. Therefore, the proof of Theorem \ref{thm:Nvar_convex} is broken down into several smaller theorems, namely Theorems \ref{thm:md_integration1}, \ref{thm:md_integration2} and \ref{thm:md_integration3}, which will be stated and proven shortly.

\begin{figure}[!ht]
\begin{tikzpicture}[node distance=1cm, auto,]
\node[punkt] (convexn) {general convex combination problem};
\node[punkt, left=2cm of convexn] (arithn) {general arithmetic mean problem} edge[pil] (convexn.west);
\node[punkt2, above=2cm of convexn] (convex2) {convex combination problem in 2 variables} edge[pil] (convexn.north);
\node[punkt, above=2cm of arithn] (arith2) {arithmetic mean problem in 2 variables} edge[pil] (convex2.west) edge[pil] (arithn.north);
\end{tikzpicture}
\caption{How our special cases of the convex combination problem relate to one another.}
\label{fig:convex_diagram}
\end{figure}
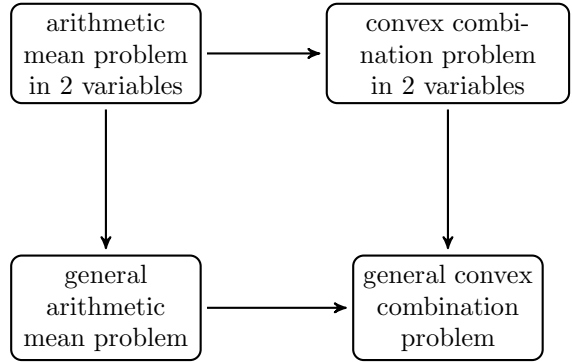

Before that, we make a short digression to review how the statement of Theorem \ref{thm:Nvar_convex} is motivated by the solutions of the convex combination problem for some special cases, which are of interest in their own right. Which special cases will be considered and how they relate to one another is shown on the diagram in Figure \ref{fig:convex_diagram}.

We start with an example that is the starting point for the arithmetic mean problem in two variables.

\begin{example}\label{ex:motivation_2var_between}
Consider the Herglotz-Nevanlinna function $q$, given by $q(z) := -\frac{1}{z}$. It is easy to check that this function is represented by the data $(0,0,\pi\delta_0)$. The functions $\til{q}_1$ and $\til{q}_2$, given by $\til{q}_1(z_1,z_2) := -\frac{1}{z_1}$ and $\til{q}_2(z_1,z_2) := -\frac{1}{z_2}$ are then represented by the data $(0,(0,0),\pi\delta_0 \otimes \lambda_\R)$ and $(0,(0,0),\lambda_\R \otimes \pi\delta_0)$, respectively, as discussed in Example \ref{ex:motivation_2var}. We note here that $\delta_0$ denotes the Dirac measure supported in the point $0 \in \R$, while $\lambda_\R$ denotes, as always, the Lebesgue measure on $\R$.

The function $\til{q}$, which is most likely to be considered as lying "halfway" between the functions $\til{q}_1$ and $\til{q}_2$, is then given by
$$\til{q}(z_1,z_2) := q(\tfrac{1}{2}z_1 + \tfrac{1}{2}z_2) = \frac{-2}{z_1+z_2}.$$
The data used to represent the function $\til{q}$ in the sense of Theorem \ref{thm:intRep_Nvar} can be verified to be $(0,0,\til{\mu})$, where the measure $\til{\mu}$ is a positive Borel measure on $\R^2$, for which the $\mu$-mass of a Borel measurable subset $U \subseteq \R^2$ is given by
$$\til{\mu}(U) := 2\pi\int_\R\chi_U(-t,t)\diff t = 2\int_\R\left(\int_\R\chi_U(t_1-t_2,t_1+t_2)\diff t_2\right)\diff(\pi\delta_0)(t_1).$$
Note that this further exemplifies the special subclass of boundary measure, considered in the previous section.\hfill$\lozenge$
\end{example}

What is perhaps most surprising about Example \ref{ex:motivation_2var_between} is that it is, in some sense, universal. 

\begin{prop}\label{prop:2var_arithm}
Let $q$ be a Herglotz-Nevanlinna function in one variable, represented by the data $(a,b,\mu)$. Then, the function $\til{q}\colon\C^{+2} \to \C$, defined by
$$\til{q}\colon (z_1,z_2) \mapsto q(\tfrac{1}{2}z_1 + \tfrac{1}{2}z_2),$$
is a Herglotz-Nevanlinna function represented by the data $(\til{a},\tilvec{b},\til{\mu})$, where $\til{a} = a$, $\tilvec{b} = (\tfrac{1}{2}b,\tfrac{1}{2}b)$ and $\til{\mu}$ is a positive Borel measure on $\R^2$, defined for any Borel measurable subset $U \subseteq \R^2$ as
\begin{equation}\label{eq:measure_arith_2var}
\til{\mu}(U) := 2\int_\R\left(\int_\R\chi_U(t_1 - t_2,t_1+t_2)\diff t_2\right)\diff\mu(t_1).
\end{equation}
\end{prop}

Regarding the proof of Proposition \ref{prop:2var_arithm}, one can show that plugging the parameters $\til{a},\tilvec{b}$ and $\til{\mu}$ into representation \eqref{eq:intRep_Nvar} gives back the function $\til{q}$, relying afterwards on the uniqueness statement of Remark \ref{rem:intRep_unique}. This is omitted as it follows directly from the more general results presented shortly. Note, however, that Theorem \ref{thm:measure_2var} guarantees that the measure $\til{\mu}$, defined by relation \eqref{eq:measure_arith_2var}, is the representing measure of some Herglotz-Nevanlinna function in two variables.

One may now continue in one of two directions, either by considering the convex combination problem in two variables, or by moving on to the general arithmetic mean problem. Choosing the former, suppose that $k_1,k_2 > 0$ such that $k_1 + k_2 = 1$. We choose now to write these coefficients as
\begin{equation}
\label{eq:convex_coeff_2var}
k_1 = \frac{1}{1+b_1}, \quad k_2 = \frac{b_1}{1+b_1}, \quad\text{or equivalently,}\quad b_1 = \frac{k_2}{k_1},
\end{equation}
where $b_1 > 0$. It is elementary to verify that this describes, in fact, a bijection between the sets $\{(k_1,k_2) \in \R^2~|~k_1,k_2 > 0,\:k_1+k_2 = 1\}$ and $\{b_1 \in \R~|~b_1 > 0\}$. 

We are now ready to modify Proposition \ref{prop:2var_arithm} in order to accommodate convex combinations. 

\begin{prop}\label{prop:2var_convex}
Let $q$ be a Herglotz-Nevanlinna function in one variable, represented by the data $(a,b,\mu)$. Let now $k_1,k_2 > 0$, such that $k_1 + k_2 = 1$. Then, the function $\til{q}\colon\C^{+2} \to \C$, defined by
$$\til{q}\colon (z_1,z_2) \mapsto q(k_1z_1 + k_2z_2),$$
is a Herglotz-Nevanlinna function represented by the data $(\til{a},\tilvec{b},\til{\mu})$, where $\til{a} = a$, $\tilvec{b} = (k_1b,k_2b)$ and $\til{\mu}$ is a positive Borel measure on $\R^2$, defined for any Borel measurable subset $U \subseteq \R^2$ as
\begin{equation}\label{eq:measure_convex_2var}
\til{\mu}(U) := \beta_2\int_\R\left(\int_\R\chi_U(t_1 - b_1t_2,t_1+t_2)\diff t_2\right)\diff\mu(t_1).
\end{equation}
Here, we have $b_1 := \frac{k_2}{k_1}$ and $\beta_2 := 1+b_1$.
\end{prop}

The proof of Proposition \ref{prop:2var_convex} follows very closely the proof of Proposition \ref{prop:2var_arithm}. One can again show that plugging the parameters $\til{a},\tilvec{b}$ and $\til{\mu}$ into representation \eqref{eq:intRep_Nvar} gives back the function $\til{q}$, which is omitted. However, as before, Theorem \ref{thm:measure_2var} guarantees that the measure $\til{\mu}$, defined by relation \eqref{eq:measure_convex_2var}, is the representing measure of some Herglotz-Nevanlinna function in two variables.

If we now wish to move on to the general arithmetic mean problem instead, the challenge becomes how to modify the definition of the measure $\til{\mu}$. The following turns out to be the right choice.

\begin{prop}\label{prop:Nvar_arith}
Let $q$ be a Herglotz-Nevanlinna function in one variable, represented by the data $(a,b,\mu)$. Let $n \geq 2$. Then, the function $\til{q}\colon\C^{+n} \to \C$, defined by
$$\til{q}\colon (z_1,z_2,\ldots,z_n) \mapsto q(\tfrac{1}{n}z_1 + \tfrac{1}{n}z_2 + \ldots + \tfrac{1}{n}z_n),$$
is a Herglotz-Nevanlinna function represented by the data $(\til{a},\tilvec{b},\til{\mu})$, where $\til{a} = a$, $\tilvec{b} = (\tfrac{1}{n}b,\tfrac{1}{n}b,\ldots,\frac{1}{n}b)$ and $\til{\mu}$ is a positive Borel measure on $\R^n$, defined for any Borel measurable subset $U \subseteq \R^n$ as
\begin{multline}\label{eq:measure_arith_Nvar}
\til{\mu}(U) \\ := n\int_\R\left(\int_{\R^{n-1}}\chi_U(t_1 - t_2,\ldots,t_1-t_n,t_1+t_2+\ldots+t_n)\diff t_n \ldots \diff t_2\right)\diff\mu(t_1).
\end{multline}
\end{prop}

Proving Proposition \ref{prop:Nvar_arith} becomes significantly more difficult, comparing with the proofs of the previous two propositions. When the parameters $\til{a},\tilvec{b}$ and $\til{\mu}$ are plugged back into representation \eqref{eq:intRep_Nvar}, we are now faced with a sequential process of $n-1$ integrations with respect to $\lambda_\R$. For these reasons, we have, so far, avoided doing any explicit calculations and will instead, as mentioned previously, present them in full only for the most general case. Moreover, we can no longer rely on Theorem \ref{thm:measure_2var}, and, as such, do not know from the beginning whether the measure $\til{\mu}$, defined using relation \eqref{eq:measure_arith_Nvar}, is the representing measure of some Herglotz-Nevanlinna function.

We have, thus, arrived at our final frontier, namely, how to combine the solutions of the arithmetic mean problem and the convex combination problem in two variables into a solution of the general convex combination problem. We begin by defining a matrix $\mat{M}_n$, for $n \geq 2$ and given numbers $b_1,b_2,\ldots,b_{n-1} > 0$, as
\begin{equation}
\label{eq:convex_matrix}
\mat{M}_n := \left[
\begin{array}{ccccc}
1 & -b_1 & ~ & ~ & ~ \\
1 & ~ & -b_2 & ~ & ~ \\
\vdots & ~ & ~ & \ddots & ~ \\
1 & ~ & ~ & ~ & -b_{n-1} \\
1 & 1 & 1 & \ldots & 1
\end{array}
\right]_{n \times n}.
\end{equation}
We note here that all the empty places in the matrix $M_n$ are filled with zeros.

This particular choice of a matrix should not be surprising, since it, for appropriate $n$ and $b_j$, describes precisely how the integration variables are intertwined in the formulas \eqref{eq:measure_arith_2var}, \eqref{eq:measure_convex_2var} and \eqref{eq:measure_arith_Nvar}. We now introduce the number $\beta_n$ as the determinant of the matrix $\mat{M}_n$, and it is an easy exercise in linear algebra to verify that
\begin{equation}
\label{eq:betaN}
\beta_n := \det(\mat{M}_n) = \sum_{j=1}^{n-1}\prod_{\substack{i = 1 \\ i \neq j}}^{n-1}b_i + \prod_{i=1}^{n-1}b_i.
\end{equation}
We note also, for example, that $\beta_2 = 1+b_1$, as it was in Proposition \ref{prop:2var_convex}.

Suppose now that $k_1,k_2,\ldots,k_n > 0$ are such that $k_1 + k_2 + \ldots + k_n = 1$. The mapping between the numbers $b_j$ and the numbers $k_\ell$ is then chosen as
\begin{equation}
\label{eq:convex_coeff_Nvar}
k_\ell = \frac{\prod_{i=1}^{n-1}b_i}{b_\ell\beta_n},\:k_n = \frac{\prod_{i=1}^{n-1}b_i}{\beta_n}, \quad\text{or equivalently,}\quad b_j = \frac{k_n}{k_j},
\end{equation}
where, $j,\ell=1,2,\ldots,n-1$. As in the case $n=2$, it is easy to check that the relations \eqref{eq:convex_coeff_Nvar} constitute a bijection between the sets
$$\{(k_1,k_2,\ldots,k_n) \in \R^n~|~k_\ell > 0,\:k_1+k_2+\ldots+k_n=1\}$$
and
$$\{(b_1,b_2,\ldots,b_{n-1}) \in \R^{n-1}~|~b_j > 0\}.$$
Therefore, if we are given the coefficients of a convex combination, we associate a positive Borel measure $\til{\mu}$ on $\R^n$ to these coefficients through the numbers $b_j$ as stated previously in formula \eqref{eq:measure_convex_Nvar}, \ie for any Borel measurable subset $U \subseteq \R^n$ we define
\begin{multline*}
\til{\mu}(U) \\ := \beta_n\int_\R\left(\int_{\R^{n-1}}\chi_U(t_1 - b_1t_2,\ldots,t_1-b_{n-1}t_n,t_1+t_2+\ldots+t_n)\diff t_n\ldots \diff t_2\right)\diff\mu(t_1).
\end{multline*}
Observe that this definitions is, of course, dependent on the underlying function $q$ of the convex combination problem, which manifests itself through its representing measure $\mu$.

Finally, we investigate what happens to the kernel $K_n$, written in the form \eqref{eq:kernel_Kn_v2}, when integrated with respect to a measure $\til{\mu}$ of the form \eqref{eq:measure_convex_Nvar}. First, we define a new kernel $\til{K}_n^0$ as
\begin{equation}
\label{eq:kernel_N0}
\til{K}_n^0(\vec{z},\vec{t}) := K_n(\vec{z},\mat{M}_n\vec{t}),
\end{equation}
where $\vec{z} \in \C^{+n}$ and $\vec{t} \in \R^n$. Here, the upper index zero is used to note that this kernel has, so far, not been integrated with respect to $\lambda_\R$ in any variable. Thus, it holds that
$$\int_{\R^n}K_n(\vec{z},\vec{t})\diff\til{\mu}(\vec{t}) = \beta_n\int_\R\left(\int_{\R^{n-1}}\til{K}_n^0(\vec{z},\vec{t})\diff t_n \ldots \diff t_2\right)\diff \mu(t_1).$$
In order to be able to explicitly evaluate the above integral, we need a general description of what happens to the kernel $\til{K}_n^0$ after it is integrated a few times with respect to the Lebesgue measure. This requires the introduction of some notation.

First, we let $m \geq 2$ and $d \geq 0$ be such that $m+d=n$. Also, let the numbers $b_j$ be as before. Consider now the factors $A_k,B_k,C_k$ and $D_k$, defined as
$$\begin{array}{LCL}
A_k := \prod_{j=2}^{k}(t_1-b_{j-1}t_j-\I), & ~ & B_k := \prod_{j=1}^k(z_j+\I), \\[0.2cm]
C_k := \prod_{j=2}^{k}(t_1-b_{j-1}t_j-z_{j-1}), & ~ & D_k := \prod_{j=2}^{k}(t_1-b_{j-1}t_j+\I).
\end{array}$$
Here, we take $\vec{z} \in \C ^{+n}$ and $\vec{t} \in \R^n$ as usual. The definitions of the factors $A_k,B_k,C_k$ and $D_k$ are valid for $1 \leq k \leq n$, while noting that empty products are equal to $1$ by convention. These factors can be thought of as the building blocks of the kernel  $K_n$ with regards to formula \eqref{eq:kernel_Kn_v2}, but transformed with respect to the definition of the kernel $\til{K}_n^0$. Furthermore, we define constants $F_m^d$ as
$$F_m^d := 1 + \sum_{j=1}^d\frac{1}{b_{m+d-j}}.$$
Constants of this particular form will appear frequently when preforming calculations using the residue theorem. Two additional expression will be useful when doing calculations, namely
$$T_m^d := F_m^d t_1 + t_2 + \ldots + t_m$$
and
$$Z_m^d := \frac{z_m}{b_m} + \frac{z_{m+1}}{b_{m+1}} + \ldots + \frac{z_{m+d-1}}{b_{m+d-1}} + z_{m+d}.$$

Finally, we introduce a notation to write down fractions, which have very long and complicated expressions as their numerator and denominator. The notations $\{\cdot/\cdot\}$ is to be understood as a fraction where anything between the symbols $\{$ and $/$ constitutes the numerator, and anything between the symbols $/$ and $\}$ constitutes the denominator. Some simple examples of the use of this notation would be
$$\big\{1\big/2\big\} = \tfrac{1}{2},\quad \big\{1+2\big/3\big\} = \tfrac{1+2}{3},\quad \big\{1\big/2+3\big\} = \tfrac{1}{2+3}.$$

We may now introduce the general kernel $\til{K}_m^d$, for $\vec{z} \in \C^{+(m+d)}$ and $\vec{t} \in \R^m$, as
\begin{multline}
\label{eq:kernel_MD_v3}
\til{K}_m^d(\vec{z},\vec{t}) := \bigg\{\I^{3m+1}A_m(T_m^d-F_m^d\I)B_{m-1}(Z_m^d+F_m^d\I) \\
-2^{m-1}F_m^d\I C_m(T_m^d-Z_m^d) \bigg/2^{m-1}A_mC_mD_m(T_m^d-F_m^d\I)(T_m^d-Z_m^d)(T_m^d+F_m^d\I)\bigg\}.
\end{multline}
Note that, when $d = 0$ and $m = n$, formula \ref{eq:kernel_MD_v3} does indeed give back the kernel $\til{K}_n^0$ as defined in formula \eqref{eq:kernel_N0}. 

While the kernel $\til{K}_m^d$ may appear long and bulky, it possesses great mathematical beauty, as the following two theorems show.

\begin{thm}\label{thm:md_integration1}
Let $n \geq 3$, let $b_1,\ldots,b_{n-1} > 0$ and let $m,d \in \N_0$ be such that $m+d = n$ with $d \geq 0$ and $m \geq 3$. Then, it holds that
\begin{equation}
\label{eq:md_integration}
\int_\R\til{K}_m^d(\vec{z},(\vec{t},t_m))\diff t_m = \frac{\pi}{b_{m-1}}\til{K}_{m-1}^{d+1}(\vec{z},\vec{t}),
\end{equation}
where the above equality holds for any $\vec{z} \in \C^{+n}$ and any $\vec{t} \in \R^{m-1}$.
\end{thm}

\proof
If we want to do any sort of calculations, then formula \eqref{eq:kernel_MD_v3} is not particularly helpful since all the variables are hidden in the building blocks of the kernel. Therefore, we write out all the terms that explicitly contain the variable $t_m$ to get
\begin{multline}
\label{eq:kernel_MD_long}
\til{K}_m^d(\vec{z},(\vec{t},t_m)) = \bigg\{\I^{3m+1}A_{m-1}(t_1-b_{m-1}t_m-\I)(F_m^dt_1+t_2+\ldots+t_m-F_m^d\I) \\
\cdot B_{m-1}(Z_m^d+F_m^d\I)-2^{m-1}F_m^d\I C_{m-1}(t_1-b_{m-1}t_m-z_{m-1}) \\
\cdot (F_m^dt_1+t_2+\ldots+t_m-Z_m^d)\bigg/2^{m-1}A_{m-1}C_{m-1}D_{m-1}(t_1-b_{m-1}t_m-\I) \\
\cdot (t_1-b_{m-1}t_m-z_{m-1})(t_1-b_{m-1}t_m+\I)(F_m^dt_1+t_2+\ldots+t_m-F_m^d\I) \\
\cdot (F_m^dt_1+t_2+\ldots+t_m-Z_m^d)(F_m^dt_1+t_2+\ldots+t_m+F_m^d\I)\bigg\}.
\end{multline}

We observe now that the kernel $\til{K}_m^d$ has six singularities with respect to the variable $t_m$, with three lying in the upper half-plane and three lying in lower half-plane. As such, we may attempt to use standard residue calculus in one complex variable in order to evaluate the left-hand side of equality \eqref{eq:md_integration}.

With respect to the variable $t_m$, the numerator of the expression $\til{K}_m^d(\vec{z},(\vec{t},t_m))$ is a polynomial of degree at most 2, while the denominator of the same expression is a polynomial of degree 6. The latter follows from observation that the leading coefficient of denominator, in this regard, is equal to
$$2^{m-1}A_{m-1}C_{m-1}D_{m-1}(-1)^3b_{m-1}^3,$$
which is non-zero due to the fact that the factors $A_{m-1},C_{m-1}$ and $D_{m-1}$ only take non-real values by definition. This shows, in particular, that the left-hand side of \eqref{eq:md_integration} is well-defined.

Let now $\Gamma_R^-$ be standard lower half-circle contour in $\C$, as was specified in the proof of Proposition \ref{prop:measure_nevan}. Due to particular rational form of the expression $\til{K}_m^d(\vec{z},(\vec{t},t_m))$ with respect to the variable $t_m$, it holds that
\begin{multline*}
\lim\limits_{R \to \infty}\int_{\Gamma_R^-}\til{K}_m^d(\vec{z},(\vec{t},t_m))\diff t_m \\
= \lim\limits_{R \to \infty}\left(\int_{-R}^R+\int_{\gamma_R^-}\right)\til{K}_m^d(\vec{z},(\vec{t},t_m))\diff t_m = \int_\R\til{K}_m^d(\vec{z},(\vec{t},t_m))\diff t_m.
\end{multline*}
On the other hand, by the residue theorem, we have that
\begin{multline*}
\int_{\Gamma_R^-}\til{K}_m^d(\vec{z},(\vec{t},t_m))\diff t_m \\
= -2\pi\I\bigg(\Res(\til{K}_m^d(\vec{z},(\vec{t},\cdot));p_1) + \Res(\til{K}_m^d(\vec{z},(\vec{t},\cdot));p_2) + \Res(\til{K}_m^d(\vec{z},(\vec{t},\cdot));p_3)\bigg)
\end{multline*}
Here, the poles of our integrand are situated at the points
$$\begin{array}{RCL}
p_1 & := & \frac{-\I+t_1}{b_{m-1}}, \\[0.35cm]
p_2 & := & \frac{-z_{m-1}+t_1}{b_{m-1}}, \\[0.35cm]
p_3 & := & -F_m^d\I-F_m^dt_1 - t_2 - \ldots - t_{m-1}.
\end{array}$$
It is important to observe that these three points lie in $\C^-$ irrespective of the particular values of $t_j \in \R$ and $b_j > 0$.

We now continue by calculating the residue at the point $p_1$, which is equal to
$$\begin{array}{LCL}
\multicolumn{3}{L}{\Res(\til{K}_m^d(\vec{z},(\vec{t},\cdot));p_1) = \lim\limits_{t_n \to p_1}\til{K}_m^d(\vec{z},(\vec{t},t_m))(t_m - \tfrac{-\I+t_1}{b_{m-1}})} \\[0.35cm]
~ & = & \tfrac{-1}{b_{m-1}}\lim\limits_{t_n \to p_1}\til{K}_m^d(\vec{z},(\vec{t},t_m))(t_1-b_{m-1}t_m-\I) \\[0.35cm]
~ & = & -\tfrac{1}{b_{m-1}}\bigg\{-2^{m-1}F_m^d\I C_{m-1} (\I-z_{m-1})(F_m^dt_1+ t_2 + \ldots + t_{m-1} + \tfrac{t_1}{b_{m-1}} \\[0.35cm]
~ & ~ & -\tfrac{\I}{b_{m-1}} - Z_m^d)\bigg/2^{m-1}A_{m-1}C_{m-1}D_{m-1}2\I(\I-z_{m-1})(F_m^dt_1 + t_2 + \ldots \\[0.35cm]
~ & ~ & + t_{m-1} + \tfrac{t_1}{b_{m-1}} - F_m^d\I - \tfrac{\I}{b_{m-1}})(F_m^dt_1 + t_2 + \ldots + t_{m-1}+\tfrac{t_1}{b_{m-1}} \\[0.35cm]
~ & ~ & -\tfrac{\I}{b_{m-1}} - Z_m^d)(F_m^d t_1 + t_2 + \ldots + t_{m-1} + \tfrac{t_1}{b_{m-1}} + F_m^d\I - \tfrac{\I}{b_{m-1}})\bigg\} \\[0.35cm]
~ & = & \tfrac{-1}{b_{m-1}^2}\bigg\{-2^{m-2}F_m^db_{m-1}C_{m-1} \bigg/ 2^{m-1}A_{m-1}C_{m-1}D_{m-1} \\[0.35cm]
~ & ~ & \cdot (T_{m-1}^{d+1}-F_{m-1}^{d+1}\I)(T_{m-1}^{d+1}+(F_m^d-\tfrac{1}{b_{m-1}})\I)\bigg\}.
\end{array}$$
When performing this calculation, there are a few thing to take note of. Firstly, the term in the numerator of expression \eqref{eq:kernel_MD_long} that starts with the constant $\I^{3m+1}$ will tend to zero as $t_m \to p_1$. In other places, the facts that $F_m^d + \frac{1}{b_{m-1}} = F_{m-1}^{d+1}$ and
$$F_m^dt_1 + t_2 + \ldots +t_{m-1} + \tfrac{t_1}{b_{m-1}} = T_{m-1}^{d+1}$$
will simplify the encountered expression greatly. Furthermore, we notice that, in the end result, the numerator and the denominator still share some factors. Given what we expect as the end result, it is inefficient to cancel out these factors now, only to be forced to multiply them back later. On the other hand, why we have chosen to factor out an extra instance of the number $\frac{1}{b_{m-1}}$ will become clear when we calculate the residue at the point $p_3$.

Moving on to the residue at the point $p_2$, we calculate that it is equal to
$$\begin{array}{LCL}
\multicolumn{3}{L}{\Res(\til{K}_m^d(\vec{z},(\vec{t},\cdot));p_2) = \lim\limits_{t_n \to p_2}\til{K}_m^d(\vec{z},(\vec{t},t_m))(t_m - \tfrac{-z_{m-1}+t_1}{b_{m-1}})} \\[0.35cm]
~ & = & \tfrac{-1}{b_{m-1}}\lim\limits_{t_n \to p_2}\til{K}_m^d(\vec{z},(\vec{t},t_m))(t_1-b_{m-1}t_m-z_{m-1}) \\[0.35cm]
~ & = & -\tfrac{1}{b_{m-1}}\bigg\{\I^{3m+1}A_{m-1}(z_{m-1}-\I)(F_m^dt_1 + t_2 + \ldots + t_{m-1} + \tfrac{t_1}{b_{m-1}}-\tfrac{z_{m-1}}{b_{m-1}} \\[0.35cm]
~ & ~ & - F_m^d\I)B_{m-1}(Z_m^d-F_m^d\I)\bigg/2^{m-1}A_{m-1}C_{m-1}D_{m-1}(z_{m-1}-\I)(z_{m-1}+ \I)  \\[0.35cm]
~ & ~ & \cdot (F_m^dt_1 + t_2 + \ldots + t_{m-1} + \tfrac{t_1}{b_{m-1}} - \tfrac{z_{m-1}}{b_{m-1}} - F_m^d\I)(F_m^dt_1 + t_2 + \ldots + t_{m-1} \\[0.35cm]
~ & ~ & + \tfrac{t_1}{b_{m-1}} - \tfrac{z_{m-1}}{b_{m-1}} - Z_m^d)(F_m^dt_1 + t_2 + \ldots + t_{m-1} + \tfrac{t_1}{b_{m-1}} - \tfrac{z_{m-1}}{b_{m-1}} + F_m^d\I)\bigg\} \\[0.35cm]
~ & = & \tfrac{-1}{b_{m-1}^2}\bigg\{\I^{3m+1}b_{m-1}A_{m-1}B_{m-2}(Z_m^d+F_m^d\I) \bigg/ 2^{m-1}A_{m-1}C_{m-1}D_{m-1} \\[0.35cm]
~ & ~ & \cdot (T_{m-1}^{d+1}-Z_{m-1}^{d+1})(T_{m-1}^{d+1}+F_m^d\I-\tfrac{z_{m-1}}{b_{m-1}})\bigg\}.
\end{array}$$
Here, the most important thing to notice is that the term in the numerator of expression \eqref{eq:kernel_MD_long} that starts with the constant $2^{m-1}$ will tend to zero as $t_m \to p_2$.

Finally, we calculate that the residue at the point $p_3$ is equal to
$$\begin{array}{LCL}
\multicolumn{3}{L}{\Res(\til{K}_m^d(\vec{z},(\vec{t},\cdot));p_3)} \\[0.35cm]
~ & = & \lim\limits_{t_n \to p_3}\til{K}_m^d(\vec{z},(\vec{t},t_m))(F_m^dt_1 + t_2 + \ldots + t_{m-1} + F_m^d\I) \\[0.35cm]
~ & ~ & \bigg\{\I^{3m+1}A_{m-1}(t_1 - b_{m-1}(-F_m^dt_1 - t_2 - \ldots - t_{m-1} - F_m^d\I) - \I)(-2F_m^d\I) \\[0.35cm]
~ & ~ & \cdot B_{m-1}(Z_m^d+F_m^d\I) - 2^{m-1}F_m^d\I C_{m-1}(t_1 - b_{m-1}(-F_m^dt_1 - t_2 - \ldots - t_{m-1} \\[0.35cm]
~ & ~ &- F_m^d\I) - z_{m-1})(-Z_m^d-F_m^d\I)\bigg/2^{m-1}A_{m-1}C_{m-1}D_{m-1}(t_1 - b_{m-1}(-F_m^dt_1 \\[0.35cm]
~ & ~ & -t_2 - \ldots - t_{m-1} - F_m^d\I) - \I)(t_1 - b_{m-1}(-F_m^dt_1 - t_2 - \ldots - t_{m-1} - F_m^d\I) \\[0.35cm]
~ & ~ & - z_{m-1})(t_1 - b_{m-1}(-F_m^dt_1 - t_2 - \ldots - t_{m-1} - F_m^d\I) + \I)(-2F_m^d\I) \\[0.35cm]
~ & ~ & \cdot(-Z_m^d-F_m^d\I)\bigg\} \\[0.35cm]
~ & = & \tfrac{-1}{b_{m-1}^2}\bigg\{\I^{3m+1}A_{m-1}(T_{m-1}^{d+1}+(F_m^d-\tfrac{1}{b_{m-1}})\I)B_{m-1} \\[0.35cm]
~ & ~ & -2^{m-2}C_{m-1}(T_{m-1}^{d+1}+F_m^d\I-\tfrac{z_{m-1}}{b_{m-1}}) \bigg/ 2^{m-1}A_{m-1}C_{m-1}D_{m-1} \\[0.35cm]
~ & ~ & \cdot (T_{m-1}^{d+1}+(F_m^d-\tfrac{1}{b_{m-1}})\I)(T_{m-1}^{d+1}+F_m^d\I-\tfrac{z_{m-1}}{b_{m-1}})(T_{m-1}^{d+1}+F_{m-1}^{d+1}\I)\bigg\}.
\end{array}$$
Here, it is important to observe that
\begin{multline*}
t_1 - b_{m-1}(-F_m^dt_1 - t_2 - \ldots - t_{m-1} - F_m^d\I) - \I \\
= b_{m-1}(\tfrac{t_1}{b_{m-1}} + F_m^dt_1 + t_2 + \ldots + t_{m-1} + F_m^d\I - \tfrac{\I}{b_{m-1}}) \\
= b_{m-1}(T_{m-1}^{d+1} + (F_m^d + \tfrac{1}{b_{m-1}})\I).
\end{multline*}
Similar simplifications are made in the cases where the last term of the starting expression is equal to $z_{m-1}$ or $+\I$. 

It thus remains to sum together the residues at the points $p_1,p_2$ and $p_3$. We begin by observing that all three residue expression have the terms $2^{m-1},A_{m-1},C_{m-1}$ and $D_{m-1}$ in their respective denominators. Furthermore, the term $T_{m-1}^{d+1}+(F_m^d-\tfrac{1}{b_{m-1}})\I$ appears in the residue expression for the points $p_1$ and $p_3$, while the term $T_{m-1}^{d+1}+F_m^d\I-\tfrac{z_{m-1}}{b_{m-1}}$ appears in the residue expression for the points $p_2$ and $p_3$. On the other hand, the expressions $T_{m-1}^{d+1}-F_{m-1}^{d+1}\I,T_{m-1}^{d+1}-Z_{m-1}^{d+1}$ and $T_{m-1}^{d+1} + F_{m-1}^{d+1}\I$ appear only in the residue expression for the point $p_1,p_2$ and $p_3$, respectively. Thus, it holds that
$$\begin{array}{LCL}
\multicolumn{3}{L}{-2\pi\I\bigg(\Res(\til{K}_m^d(\vec{z},(\vec{t},\cdot));p_1) + \Res(\til{K}_m^d(\vec{z},(\vec{t},\cdot));p_2) + \Res(\til{K}_m^d(\vec{z},(\vec{t},\cdot));p_3)\bigg)} \\[0.35cm]
~ & = & \tfrac{2\pi\I}{b_{m-1}^2}\bigg\{(*1)\bigg/2^{m-1}A_{m-1}C_{m-1}D_{m-1}(T_{m-1}^{d+1}-F_{m-1}^{d+1}\I)  \\[0.35cm]
~ & ~ & \cdot (T_{m-1}^{d+1}+(F_m^d-\tfrac{1}{b_{m-1}})\I)(T_{m-1}^{d+1}-Z_{m-1}^{d+1})(T_{m-1}^{d+1}+F_m^d\I-\tfrac{z_{m-1}}{b_{m-1}}) \\[0.35cm]
~ & ~ & \cdot (T_{m-1}^{d+1} + F_{m-1}^{d+1}\I)\bigg\} = (**1),
\end{array}$$
where the expression $(*1)$ is given by
$$\begin{array}{LCL}
(*1) & := & -2^{m-2}F_m^db_{m-1}C_{m-1}(T_{m-1}^{d+1}-Z_{m-1}^{d+1})(T_{m-1}^{d+1}+F_m^d\I-\tfrac{z_{m-1}}{b_{m-1}}) \\[0.35cm]
~ & ~ & \cdot(T_{m-1}^{d+1}+F_{m-1}^{d+1}\I) + \I^{3m+1}b_{m-1}A_{m-1}B_{m-2}(Z_m^d+F_m^d\I)(T_{m-1}^{d+1}-F_{m-1}^{d+1}) \\[0.35cm]
~ & ~ & \cdot(T_{m-1}^{d+1}+(F_m^d-\tfrac{1}{b_{m-1}})\I)(T_{m-1}^{d+1}+F_{m-1}^{d+1}\I)+\I^{3m+1}A_{m-1}B_{m-1} \\[0.35cm]
~ & ~ & \cdot(T_{m-1}^{d+1}+(F_m^d-\tfrac{1}{b_{m-1}})\I)(T_{m-1}^{d+1}-F_{m-1}^{d+1}\I)(T_{m-1}^{d+1}-Z_{m-1}^{d+1}) \\[0.35cm]
~ & ~ & -2^{m-2}C_{m-1}(T_{m-1}^{d+1}+F_{m}^{d}\I-\tfrac{z_{m-1}}{b_{m-1}})(T_{m-1}^{d+1}-F_{m-1}^{d+1}\I)(T_{m-1}^{d+1}-Z_{m-1}^{d+1}).
\end{array}$$
Summing together the two terms of expression $(*1)$ that begin with the constant $\I^{3m+1}$ gives
$$\begin{array}{LCL}
\multicolumn{3}{L}{\I^{3m+1}A_{m-1}B_{m-2}(T_{m-1}^{d+1}-F_{m-1}^{d+1\I})(T_{m-1}^{d+1}+(F_m^d-\tfrac{1}{b_{m-1}})\I)\bigg(b_{m-1}(Z_m^d+F_m^d\I)} \\[0.35cm]
~ & ~ & \cdot(T_{m-1}^{d+1}+F_{m-1}^{d+1}\I)+(z_{m-1}+\I)(T_{m-1}^{d+1}-Z_{m-1}^{d+1}\I)\bigg) \\[0.35cm]
~ & = & \I^{3m+1}A_{m-1}B_{m-2}(T_{m-1}^{d+1}-F_{m-1}^{d+1}\I)(T_{m-1}^{d+1}+(F_m^d-\tfrac{1}{b_{m-1}})\I)b_{m-1} \\[0.35cm]
~ & ~ & \cdot(T_{m-1}^{d+1}+F_m^d\I-\tfrac{z_{m-1}}{b_{m-1}})(Z_{m-1}^{d+1}+F_{m-1}^{d+1}\I),
\end{array}$$
while summing together the two terms of expression $(*1)$ that begin with the constant $2^{m-2}$ gives
$$\begin{array}{LCL}
\multicolumn{3}{L}{-2^{m-2}C_{m-1}(T_{m-1}^{d+1}-Z_{m-1}^{d+1})(T_{m-1}^{d+1}+F_{m}^{d}\I-\tfrac{z_{m-1}}{b_{m-1}})\bigg(b_{m-1}F_m^d} \\[0.35cm]
~ & ~ & \cdot(T_{m-1}^{d+1}+F_{m-1}^{d+1}\I)+(T_{m-1}^{d+1}-F_{m-1}^{d+1}\I)\bigg) \\[0.35cm]
~ & = & -2^{m-2}C_{m-1}(T_{m-1}^{d+1}-Z_{m-1}^{d+1})(T_{m-1}^{d+1}+F_{m}^{d}\I-\tfrac{z_{m-1}}{b_{m-1}})b_{m-1}F_{m-1}^{d+1} \\[0.35cm]
~ & ~ & \cdot(T_{m-1}^{d+1}+(F_{m}^{d}-\tfrac{1}{b_{m-1}})\I).
\end{array}$$
Using these simplifications, we conclude that
$$\begin{array}{LCL}
(**1) & = & \tfrac{2\pi\I}{b_{m-1}^2}\bigg\{\I^{3m+1}A_{m-1}B_{m-2}b_{m-1}(T_{m-1}^{d+1}-F_{m-1}^{d+1}\I)(Z_{m-1}^{d+1}+F_{m-1}^{d+1}\I) \\[0.35cm]
~ & ~ & -2^{m-2}b_{m-1}F_{m-1}^{d+1}C_{m-1}(T_{m-1}^{d+1}-Z_{m-1}^{d+1})\bigg/2^{m-1}A_{m-1}C_{m-1}D_{m-1} \\[0.35cm]
~ & ~ & \cdot(T_{m-1}^{d+1}-F_{m-1}^{d+1}\I)(T_{m-1}^{d+1}-Z_{m-1}^{d+1})(T_{m-1}^{d+1} + F_{m-1}^{d+1}\I)\bigg\} \\[0.35cm]
~ & = & \tfrac{\pi}{b_{m-1}}\bigg\{\I^{3m+2}A_{m-1}B_{m-2}(T_{m-1}^{d+1}-F_{m-1}^{d+1}\I)(Z_{m-1}^{d+1}+F_{m-1}^{d+1}\I) \\[0.35cm]
~ & ~ & -2^{m-2}F_{m-1}^{d+1}\I C_{m-1}(T_{m-1}^{d+1}-Z_{m-1}^{d+1})\bigg/2^{m-2}A_{m-1}C_{m-1}D_{m-1} \\[0.35cm]
~ & ~ & \cdot(T_{m-1}^{d+1}-F_{m-1}^{d+1}\I)(T_{m-1}^{d+1}-Z_{m-1}^{d+1})(T_{m-1}^{d+1} + F_{m-1}^{d+1}\I)\bigg\} \\
~ & = & \tfrac{\pi}{b_{m-1}}\til{K}_{m-1}^{d+1}(\vec{z},\vec{t}).
\end{array}$$
Here, the last equality follows from the observation that $\I^{3m+2} = \I^{3m-2}$. This finishes the proof.
\endproof

\begin{thm}\label{thm:md_integration2}
Let $n \geq 2$ and let $b_1,\ldots,b_{n-1} > 0$. Then, it holds that
\begin{equation}
\label{eq:md_integration_final}
\int_\R\til{K}_2^{n-2}(\vec{z},(t_1,t_2))\diff t_2 = \pi\frac{\prod_{j=2}^{n-1}b_j}{\beta_n}K_1(k_1z_1+k_2z_2+\ldots+k_nz_n,t_1),
\end{equation}
where the above equality holds for any $\vec{z} \in \C^{+n}$ and any $t_1 \in \R$. Here, the numbers $b_j$ and $k_\ell$ are related by formula \eqref{eq:convex_coeff_Nvar}.
\end{thm}

\proof
In short, the proof of this theorem is exactly the same as the proof of the preceding theorem. More precisely, all but the very last calculations, performed in the proof of Theorem \ref{thm:md_integration1}, turn out to still be valid, even now when $m=2$ and $d = n-2$. In particular, the same arguments as before justify the use of residue theorem, and the residues at the points $p_1,p_2$ and $p_3$ are still given by the same expressions as before, implying that
$$\begin{array}{LCL}
\multicolumn{3}{L}{\int_\R\til{K}_2^{n-2}(\vec{z},(t_1,t_2))\diff t_2 = \lim\limits_{R \to \infty}\int_{\Gamma_R^-}\til{K}_2^{n-2}(\vec{z},(t_1,t_2))\diff t_2} \\[0.35cm]
~ & = & -2\pi\I\bigg(\Res(\til{K}_2^{n-2}(\vec{z},(t_1,\cdot));\tfrac{-\I+t_1}{b_1}) + \Res(\til{K}_2^{n-2}(\vec{z},(t_1,\cdot));\tfrac{-z_1+t_1}{b_1}) \\[0.35cm]
~ & ~ & + \Res(\til{K}_2^{n-2}(\vec{z},(t_1,\cdot));-F_2^{n-2}t_1 - F_2^{n-2}\I)\bigg) \\[0.35cm]
~ & = & \tfrac{2\pi\I}{b_{1}^2}\bigg\{(*2)\bigg/2A_{1}C_{1}D_{1}(T_{1}^{n-1}-F_{1}^{n-1}\I)(T_{1}^{n-1}+(F_2^{n-2}-\tfrac{1}{b_{1}})\I)  \\[0.35cm]
~ & ~ & \cdot (T_{1}^{n-1}-Z_{1}^{n-1})(T_{1}^{n-1}+F_2^{n-2}\I-\tfrac{z_{1}}{b_{1}})(T_{1}^{n-1} + F_{1}^{n-1}\I)\bigg\} = (**2),
\end{array}$$
where the expression $(*2)$ is, after analogous simplifications as before, given by
$$\begin{array}{LCL}
(*2) & := & -\I b_1(T_1^{n-1}-F_1^{n-1}\I)(T_1^{n-1}+(F_2^{n-1}-\tfrac{1}{b_1})\I)(T_1^{n-1}+F_2^{n-2}\I-\tfrac{z_1}{b_1}) \\[0.35cm]
~ & ~  & \cdot(Z_1^{n-1}+F_1^{n-1}\I)-b_1F_1^{n-1}(T_1^{n-1}-Z_1^{n-1})(T_1^{n-1}+(F_2^{n-1}-\tfrac{1}{b_1})\I) \\[0.35cm]
~ & ~ & \cdot(T_1^{n-1}+F_2^{n-2}\I-\tfrac{z_1}{b_1}).
\end{array}$$
However, unlike in the proof of Theorem \ref{thm:md_integration1}, these expressions can be further simplified using the observations that $A_1 = C_1 = D_1 = 1$ and that
$$T_1^{n-1} \pm F_1^{n-1}\I = F_1^{n-1}t_1 \pm F_1^{n-1}\I = F_1^{n-1}(t_1 \pm \I),$$
as well as $T_1 - Z_1^{n-1} = F_1^{n-1}t_1 - Z_1^{n-1}$. Thus, it holds that
$$\begin{array}{LCL}
\multicolumn{3}{L}{(**2) = \tfrac{2\pi\I}{b_1^2}\bigg\{-\I b_1 F_1^{n-1}(t_1 - \I)(Z_1^{n-1}+F_1^{n-1}\I)-b_1F_1^{n-1}(F_1^{n-1}t_1 - Z_1^{n-1})\bigg/} \\[0.35cm]
~ & ~ & 2(F_1^{n-1})^2(t_1-\I)(t_1+\I)(F_1^{n-1}t_1 - Z_1^{n-1})\bigg\} \\[0.35cm]
~ & = & \tfrac{\pi}{b_1F_1^{n-1}}\bigg\{(t_1 - \I)(Z_1^{n-1}+F_1^{n-1}\I)-\I(F_1^{n-1}t_1-Z_1^{n-1})\bigg/ \\[0.35cm]
~ & ~ & (1+t_1^2)(F_1^{n-1}t_1-Z_1^{n-1})\bigg\} \\[0.35cm]
~ & = & \tfrac{\pi}{b_1F_1^{n-1}}\bigg\{(F_1^{n-1}+t_1Z_1^{n-1})\bigg/(1+t_1^2)(F_1^{n-1}t_1-Z_1^{n-1})\bigg\} \\[0.35cm]
~ & = & \tfrac{\pi}{b_1F_1^{n-1}}\bigg\{(1+t_1\tfrac{Z_1^{n-1}}{F_1^{n-1}})\bigg/(1+t_1^2)(t_1-\tfrac{Z_1^{n-1}}{F_1^{n-1}})\bigg\} = \tfrac{\pi}{b_1F_1^{n-1}}K_1(\tfrac{Z_1^{n-1}}{F_1^{n-1}},t_1).
\end{array}$$
Observing that
$$\frac{1}{b_1F_1^{n-1}} = \frac{\prod_{j=2}^{n-1}b_j}{\beta_n} \quad\text{and}\quad \frac{Z_1^{n-1}}{F_1^{n-1}} = k_1z_1 + k_2z_2 + \ldots +k_nz_n$$
finishes the proof.
\endproof

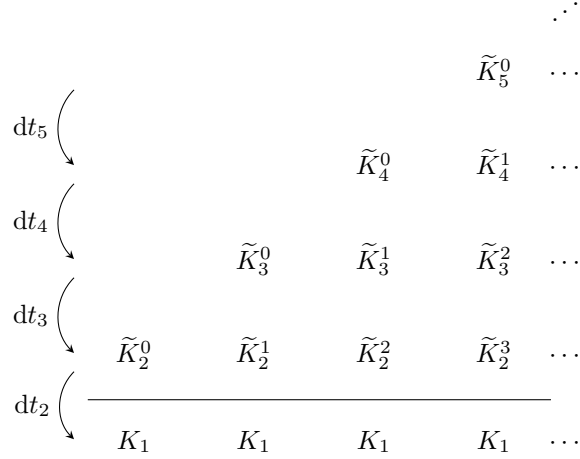
\begin{figure}[!ht]
\begin{tikzpicture}
\matrix (m) [matrix of math nodes,row sep=0.5em,column sep=0.75em,minimum width=0.25em]{
~ & ~ & ~ & ~ & ~ & ~ & ~ & ~ & \iddots \\
~ & ~ & ~ & ~ & ~ & ~ & ~ & \til{K}_5^0  & \cdots \\
~ & ~ & ~ & ~ & ~ & ~ & ~ & ~  & ~  \\
~ & ~ & ~ & ~ & ~ & \til{K}_4^0 & ~ & \til{K}_4^1  & \cdots \\
~ & ~ & ~ & ~ & ~ & ~ & ~ & ~  & ~  \\
~ & ~ & ~ & \til{K}_3^0 & ~ & \til{K}_3^1 & ~ & \til{K}_3^2  & \cdots \\
~ & ~ & ~ & ~ & ~ & ~ & ~ & ~  & ~  \\
~ & \til{K}_2^0 & ~ & \til{K}_2^1 & ~ & \til{K}_2^2 & ~ & \til{K}_2^3  & \cdots \\
~ & ~ & ~ & ~ & ~ & ~ & ~ & ~  & ~  \\
~ & K_1 & ~ & K_1 & ~ & K_1 & ~ & K_1  & \cdots \\
};
\path[-stealth]
(m-9-1.east) edge [solid,-] (m-9-9)
(m-8-1.south) edge [out=225,in=135] node [left] {$\diff t_2$} (m-10-1.north)
(m-6-1.south) edge [out=225,in=135] node [left] {$\diff t_3$} (m-8-1.north)
(m-4-1.south) edge [out=225,in=135] node [left] {$\diff t_4$} (m-6-1.north)
(m-2-1.south) edge [out=225,in=135] node [left] {$\diff t_5$} (m-4-1.north);
\end{tikzpicture}
\caption{The kernel ladder.}
\label{fig:kernel_ladder}
\end{figure}

The statements of Theorems \ref{thm:md_integration1} and \ref{thm:md_integration2} can thought of as a, sort of, ladder, visualized in Figure \ref{fig:kernel_ladder}. Each column, depicted in Figure \ref{fig:kernel_ladder}, shows the number of integrations of the kernel $\til{K}_m^d$ with respect to the Lebesgue measure, described, in total, by Theorems \ref{thm:md_integration1} and \ref{thm:md_integration2}. If $n=2$, then there is no need for Theorem \ref{thm:md_integration1}, and we only need to apply Theorem \ref{thm:md_integration2}. If $n=3$, we can first apply Theorem \ref{thm:md_integration1} once, followed by Theorem \ref{thm:md_integration2}. In general, the reduction of the kernel $\til{K}_n^0$ to the kernel $K_1$ is summarized by the following theorem.

\begin{thm}\label{thm:md_integration3}
Let $n \geq 2$ and let $b_1,\ldots,b_{n-1} > 0$. Then, it holds that
\begin{equation}
\label{eq:md_integration_complete}
\int_{\R^{n-1}}\til{K}_n^0(\vec{z},(t_1,t_2,\ldots,t_n))\diff t_n \ldots \diff t_2 = \frac{\pi^{n-1}}{\beta_n}K_1(k_1z_1+k_2z_2+\ldots+k_nz_n,t_1),
\end{equation}
where the above equality holds for any $\vec{z} \in \C^{+n}$ and any $t_1 \in \R$. Here, the numbers $b_j$ and $k_\ell$ are related by formula \eqref{eq:convex_coeff_Nvar}.
\end{thm}

\proof
First, we apply Theorem \ref{thm:md_integration1} sequentially $(n-2)$-times on the left-hand side of equality \eqref{eq:md_integration_complete}. Afterwards, we apply Theorem \ref{thm:md_integration2} to arrive at the end result.
\endproof

We are now ready to prove that Theorem \ref{thm:Nvar_convex} solves the convex combination problem.

\proof
\textit{(of Theorem \ref{thm:Nvar_convex}.)} We recall first what we want to prove. We have $n \geq 2$ and numbers $k_\ell > 0$ for $\ell=1,2,\ldots,n$, such that $k_1 + k_2 + \ldots + k_n = 1$.  We also have a Herglotz-Nevanlinna function $q$, represented by the data $(a,b,\mu)$, and we have defined the function $\til{q}$ by setting
$$\til{q}\colon (z_1,z_2,\ldots,z_n) \mapsto q(k_1z_1 + k_2z_2 + \ldots +k_nz_n).$$
It is trivial to see $\til{q}$ is also a Herglotz-Nevanlinna function. As such, it is represented by some data $(\til{a},\tilvec{b},\til{\mu})$ in the sense of Theorem \ref{thm:intRep_Nvar}. We claim that $\til{a} = a$, $\tilvec{b} = (k_1b,k_2b,\ldots,k_nb)$ and $\til{\mu}$ is a positive Borel measure on $\R^n$, given by formula \eqref{eq:measure_convex_Nvar}. The numbers $b_j$ and $k_\ell$ are related, as usual, by formulas \eqref{eq:convex_coeff_Nvar}, while the number $\beta_n$ was defined by relation \eqref{eq:betaN}.

Let us see now what happens if we try to integrate the kernel $K_n$ with respect to the measure $\til{\mu}$. Since we do not know is the measure $\til{\mu}$ satisfies the growth condition \eqref{eq:growth_Nvar} or the Nevanlinna condition \eqref{eq:Nevan_Nvar}, we do not know if the result will be a Herglotz-Nevanlinna function. In fact, we do not even know if the integral of the kernel $K_n$ with respect to the measure $\til{\mu}$ is well-defined. Nevertheless, we may try and see what happens, yielding first that
$$\int_{\R^n}K_n(\vec{z},\vec{t})\diff\til{\mu}(\vec{t}) = \beta_n\int_\R\left(\int_{\R^{n-1}}\til{K}_n^0(\vec{z},\vec{t})\diff t_n \ldots \diff t_2\right)\diff \mu(t_1) = (*)$$
by the definition of the measure $\til{\mu}$. Here, we note also that the definition of the measure $\til{\mu}$ is such that we first integrate the kernel $\til{K}_n^0$ with respect to $\diff t_n$, followed by $\diff t_{n-1}$, and so forth. But each of these integrals is, in fact, well-defined by Theorems \ref{thm:md_integration1} and \ref{thm:md_integration2}, with the final result, by Theorem \ref{thm:md_integration3}, being equal to
$$(*) = \pi^{n-1}\int_\R K_1(k_1z_1 + k_2z_2 + \ldots + k_nz_n,t_1)\diff \mu(t_1).$$

We may now plug in the parameters $\til{a},\tilvec{b}$ and $\til{\mu}$ into representation \eqref{eq:intRep_Nvar}, yielding
\begin{multline*}
a + k_1bz_1 + k_2bz_2 + \ldots + k_nbz_n + \frac{1}{\pi^n}\int_{\R^n}K_n(\vec{z},\vec{t})\diff\til{\mu}(\vec{t}) \\[0.35cm]
= a + b(k_1z_1 + k_2z_2 + \ldots + k_nz_n) + \frac{1}{\pi}\int_\R K_1(k_1z_1 + k_2z_2 + \ldots + k_nz_n,t_1)\diff \mu(t_1) \\[0.35cm]
= q(k_1z_1 + k_2z_2 + \ldots + k_nz_n) = \til{q}(z_1,z_2,\ldots,z_n).
\end{multline*}
Thus, the parameters $\til{a},\tilvec{b}$ and $\til{\mu}$ yield a Herglotz-Nevanlinna function when inserted into representation \eqref{eq:intRep_Nvar} and, by the uniqueness Remark \ref{rem:intRep_unique}, they must be equal to the data of the function $\til{q}$. This finishes the proof.
\endproof

We present now two corollaries that follow immediately from the proof of Theorem \ref{thm:Nvar_convex}.

\begin{coro}
Any positive Borel measure $\til{\mu}$ on $\R^n$ of the form \eqref{eq:measure_convex_Nvar} for some numbers $b_j > 0,\: j=1,2,\ldots,n-1$, satisfies both the growth condition \eqref{eq:growth_Nvar} and the Nevanlinna condition \eqref{eq:Nevan_Nvar}.
\end{coro}

\begin{coro}\label{coro:intRep_convex}
A Herglotz-Nevanlinna function $\til{q}$, constructed from a Herglotz-Nevanlinna function $q$ as in Theorem \ref{thm:Nvar_convex}, admits an integral representation formula of the from
\begin{multline}
\label{eq:intRep_convex}
\til{q}(z_1,z_2,\ldots,z_n) \\[0.35cm]
= a + k_1bz_1 + k_2bz_2 + \ldots + k_nbz_n + \frac{\beta_n}{\pi^n}\int_\R\left(\int_{\R^{n-1}}\til{K}_n^0(\vec{z},\vec{t})\diff t_n \ldots \diff t_2\right)\diff \mu(t_1),
\end{multline}
where all the parameters are as in Theorem \ref{thm:Nvar_convex}.
\end{coro}

Finally, we remove the little constraint that has been with us since the beginning of this section, namely that all the numbers $k_\ell$ have to be positive.

\begin{coro}\label{coro:intRep_convex_general}
Let $q$ be a Herglotz-Nevanlinna function in one variable, represented by the data $(a,b,\mu)$. Let now $n \geq 2$ and $k_\ell \geq 0$ for $\ell=1,2,\ldots,n$, such that $k_1 + k_2 + \ldots +k_n = 1$. Let
$$R := \{i_1,i_2,\ldots,i_\rho\} = \{i_\ell~|~k_{i_\ell} = 0\} \subseteq \{1,2,\ldots,n\}$$ 
and
$$S := \{j_1,j_2,\ldots,j_\sigma\} = \{1,2,\ldots,n\}\setminus R$$
be sets of sizes $\rho$ and $\sigma$, respectively, where we assume that the elements of these sets are indexed in ascending order. Then, the function $\hat{q}\colon \C^{+n} \to \C$, defined by
$$\hat{q}\colon (z_1,z_2,\ldots,z_n) \mapsto q(k_1z_1 + k_2z_2 + \ldots +k_nz_n),$$
is a Herglotz-Nevanlinna function represented by the data $(\hat{a},\hatvec{b},\hat{\mu})$, where $\hat{a} = a$, $\hatvec{b} = (k_1b,k_2b,\ldots,k_nb)$ and $\hat{\mu}$ is a positive Borel measure on $\R^n$, defined for any Borel measurable subset $U \subseteq \R^n$ as
\begin{equation*}\label{eq:measure_convex_general_Nvar}
\hat{\mu}(U) := \int_{\R^\sigma}\left(\int_{\R^\rho}\chi_U(\vec{t})\diff t_{i_\rho}\diff t_{i_{\rho-1}} \ldots \diff t_{i_1}\right)\diff\til{\mu}((t_{j_1},t_{j_2},\ldots,t_{j_\sigma})).
\end{equation*}
Here, the measure $\til{\mu}$ is taken as the representing measure of the Herglotz-Nevanlinna function
$$\til{q}\colon (z_1,z_2,\ldots,z_\sigma) \mapsto q(k_{j_1}z_1 + k_{j_2}z_2 + \ldots + k_{j_\sigma}z_\sigma).$$
\end{coro}

In short, the above corollary states that, if we allow some coefficients $k_\ell$ to be zero, we should first integrate out the corresponding $t$-variables and then use Theorem \ref{thm:Nvar_convex}.

\proof
We recall the fact, mentioned in Example \ref{ex:motivation_2var}, that integrating the kernel $K_n$ once with respect to $\diff t_j$ gives a constant multiple of $K_{n-1}$ with the $j$-th variable missing \cite[Example 3.4]{LugerNedic2017}. More precisely, it holds that
\begin{multline*}
\int_\R K_n((z_1,\ldots,z_n),(t_1,\ldots,t_{j-1},t_j,t_{j+1},\ldots,t_n))\diff t_j \\
= \pi K_{n-1}((z_1,\ldots,z_{j-1},z_{j+1},\ldots,z_n),(t_1,\ldots,t_{j-1},t_{j+1},\ldots,t_n)).
\end{multline*}
If $\sigma = 1$, then $\til{q}$ is a function of one variable and the the measure $\til{\mu}$ is just the measure $\mu$. Otherwise, the measure $\til{\mu}$ is as described by Theorem \ref{thm:Nvar_convex}. The result then follows.
\endproof

\section*{Acknowledgements}

The author would like to thank Alexandru Aleman for initiating this work by inquiring about the arithmetic mean problem in two variables, as well as Daniel Sj\"{o}berg for his seminary question that led to the consideration of the convex combination problem. The author would also like to thank Lars Jonsson and Annemarie Luger for enthusiastic discussions on the subject and careful reading of the manuscript.

\bibliographystyle{amsplain}

\begin{thebibliography}{99}

\bibitem{AbramowitzStegun1964}
M. Abramowitz and I. A. Stegun, \emph{Handbook of Mathematical Functions: With Formulas, Graphs, and Mathematical Tables}, Applied mathematics series, Dover Publications, 1964.

\bibitem{agler2}
J. Agler, J. E. McCarthy, N. J. Young, \emph{Operator monotone functions and L\"{o}wner functions of several variables}, Ann. of Math. (2) \textbf{176} (2012), no. 3, 1783--1826. 

\bibitem{agler}
J. Agler, R. Tully-Doyle, N. J. Young, \emph{Nevanlinna representations in several variables},  J. Funct. Anal. \textbf{270} (2016), no. 8, 3000--3046. 

\bibitem{Akhiezer1965}
N. I. Akhiezer, \emph{The classical moment problem and some related questions in analysis}, Hafner Publishing Co., New York, 1965.

\bibitem{Bernland2011}
A. Bernland, A. Luger, M. Gustafsson, \emph{Sum rules and constraints on passive systems}, J. Phys. A: Math. Theor. \textbf{44} (2011), 145205.

\bibitem{Brune1931}
O. Brune, \emph{Synthesis of a Finite Two-terminal Network whose Driving-point Impedance is a Prescribed Function of Frequency}, MIT J. Math. Phys. \textbf{10}, (1931), 191--236.

\bibitem{Cauer1932}
W. Cauer, \emph{The Poisson integral for functions with positive real part}, Bull. Amer. Math. Soc. \textbf{38} (1932), no. 10, 713--717.

\bibitem{Vladimirov1974}
Yu. N. Drozhzhinov and V. S. Vladimirov, \emph{Holomorphic functions in a polydisc with non-negative imaginary part}, Mat. Zametki \textbf{15} (1974), 55--61. English transl. in Math. Notes \textbf{15} (1974), 31--34.

\bibitem{GoldenPapanicolaou1983}
K. Golden and G. Papanicolaou, \emph{Bounds for effective parameters of heterogeneous media by analytic continuation}, Comm. Math. Phys. \textbf{90} (1983), no. 4, 473--491.

\bibitem{GoldenPapanicolaou1985}
K. Golden and G. Papanicolaou, \emph{Bounds for effective parameters of multicomponent media by analytic continuation}, J. Statist. Phys. \textbf{40} (1985), no. 5-6, 655--667.

\bibitem{Ivanenko+etal2017}
Y. Ivanenko et. al., \emph{Passive approximation and optimization with {B}-splines}, URN: urn:nbn:se:lnu:diva-63878, Department of Physics and Electrical Engineering, 351 95 V\"{a}xj\"{o}, Sweden, 2017.

\bibitem{KacKrein1974}
I. S. Kac and M. G. Krein, \emph{R-functions--analytic functions mapping the upper half-plane into itself}, Amer. Math. Soc. Transl. (2) \textbf{103} (1974), 1--18.

\bibitem{King2009}
F. W. King, \emph{Hilbert transforms vol. I-II}, Cambridge University Press, 2009.

\bibitem{LugerNedic2016}
A. Luger and M. Nedic, \emph{A characterization of Herglotz-Nevanlinna function in two variables via integral representations}, Ark. Mat. \textbf{55} (2017), 199--216.

\bibitem{LugerNedic2017}
A. Luger and M. Nedic, \emph{An integral representation for Herglotz-Nevanlinna functions in several variables}, arXiv: 1705.10562.

\bibitem{Milton2016}
G. W. Milton (editor), \emph{Extending the Theory of Composites to Other Areas of Science}, BookBaby Print, 2016.

\bibitem{Milton2002}
G. W. Milton, \emph{The Theory of Composites}, Cambridge Monographs on Applied and Computational Mathematics, Cambridge University Press, 2002.

\bibitem{Nevanlinna1922}
R. Nevanlinna, \emph{Asymptotische Entwicklungen beschr\"{a}nkter Funktionen und das Stieltjessche Momentenproblem}, Ann. Acad. Sci. Fenn. (A) \textbf{18 (5)} (1922), 1--53.

\bibitem{Simon1998}
B. Simon, \emph{The classical moment problem as a self-adjoint finite difference operator}, Adv. Math. \textbf{137} (1998), no. 1, 82--203.

\bibitem{Vladimirov1979}
V. S. Vladimirov, \emph{Generalized function in mathematical physics},  "Nauka", 1979; English transl., Mir publishers, Moscow, 1979.

\bibitem{Vladimirov1969}
V. S. Vladimirov, \emph{Holomorphic functions with non-negative imaginary part in a tubular domain over a cone},  Mat. Sb. \textbf{79 (121)} (1969), 182--152. English transl. in Math. USSR-Sb. \textbf{8} (1969), 125--146.

\end{thebibliography}

\end{document}